\documentclass[preprint,12pt]{elsarticle}
\usepackage[utf8]{inputenc}
\usepackage{textcomp}
\usepackage{pstricks}
\usepackage{amssymb}
\usepackage{pifont}
\usepackage{graphicx}
\usepackage[version=4]{mhchem}
\usepackage{siunitx}
\usepackage{longtable,tabularx}
\usepackage[numbers]{natbib}
\usepackage{hyperref}
\usepackage[capitalise]{cleveref}
\usepackage{amsfonts,amsmath,mathtools,amssymb}
\usepackage{setspace}
\usepackage{placeins}
\usepackage[colorinlistoftodos]{todonotes}
\usepackage{booktabs}
\usepackage{float}
\usepackage{algpseudocode}
\usepackage{algorithm}
\usepackage{subfigure}
\usepackage{lineno}
\usepackage{tikz}
\usepackage{multirow}
\usepackage{rotating}
\usepackage{appendix}
\usepackage{makecell}
\usepackage[makeroom]{cancel}
\usetikzlibrary{arrows,calc,through,backgrounds,matrix,decorations.pathmorphing}
\usepackage{soul}
\usepackage{rotating}
\sethlcolor{yellow}
\soulregister\ref7
\soulregister\cref7
\soulregister\cite7
\soulregister\citet7
\soulregister\Cref7
\usepackage{multirow}
\newcommand{\ra}[1]{\renewcommand{\arraystretch}{#1}}

\newcommand{\f}[2]{\frac{#1}{#2}}
\newcommand{\mb}[1]{\mathbf{#1}}

\newcommand{\mbg}[1]{\boldsymbol{\mathbf{#1}}}

\DeclareMathAlphabet\mathbfcal{OMS}{cmsy}{b}{n}
\renewcommand{\d}{\mathop{}\!\mathrm{d}} 
\newcommand{\p}{\partial}
\crefname{equation}{eq.}{eqs.}  
\Crefname{equation}{Eq.}{Eqs.}  

\Crefname{appendix}{}{}

\usepackage[dvipsnames,table,xcdraw]{xcolor}

\definecolor{rev1}{HTML}{FF999A} 
\definecolor{rev2}{HTML}{F3F298} 
\definecolor{rev3}{HTML}{B2E0AE} 
\definecolor{other}{HTML}{C8C7FF} 

\begin{document}

\title{Modal-Centric Field Inversion via Differentiable Proper Orthogonal Decomposition}

\author[inst1]{Rohit Sunil Kanchi}
\author[inst1]{Sicheng He}

\affiliation[inst1]{organization={Department of Mechanical and Aerospace Engineering,\\ University of Tennessee},
city={Knoxville},
postcode={TN 37996},
country={USA}}



\begin{frontmatter}
\begin{abstract}
\tolerance=1000
Inverse problems in computational physics often require matching high-dimensional spatio-temporal fields, leading to prohibitive computational costs and ill-conditioned optimizations.
We introduce modal-centric field inversion (MCFI), a paradigm that reformulates inverse problems in the reduced space of proper orthogonal decomposition (POD) modes rather than the full physical state space.
By targeting dominant flow structures instead of point-wise field values, MCFI provides a compact, physically meaningful objective that naturally regularizes the inversion and dramatically reduces computational burden.
Central to this framework is the differentiable POD\@: an adjoint-based method that efficiently computes sensitivities of POD modes with respect to model parameters, enabling gradient-based optimization in the modal space.
We demonstrate MCFI on a one and two-dimensional modified viscous Burger's equation, optimizing spatially varying coefficients to match target dynamics through mode-matching.
The adjoint formulation achieves computational cost independent of parameter dimension, in contrast to finite-difference approaches that scale linearly.
MCFI establishes a foundation for scalable inverse design and model calibration in unsteady, high-dimensional systems.
\end{abstract}
\end{frontmatter}

\section{Introduction}\label{sec:intro}

Inverse problems arise throughout computational science and engineering whenever model parameters, boundary conditions, or constitutive fields must be inferred from indirect observations.
Applications span heat transfer coefficient estimation from temperature measurements~\citep{beck1985inverse}, material property identification in structural mechanics~\citep{bui1994inverse}, subsurface characterization in geophysics~\citep{tarantola2005inverse}, and turbulence model calibration in fluid dynamics.
These problems share common computational challenges: high-dimensional parameter spaces, expensive forward model evaluations, ill-posedness requiring regularization, and the need for rigorous uncertainty quantification.
Solution strategies include Bayesian inference, adjoint-based optimization, ensemble methods, and increasingly, data-driven surrogates.
In the context of computational fluid dynamics (CFD), turbulence model calibration exemplifies these challenges and has motivated significant methodological advances.

Field inversion machine learning (FIML) reframes calibration as a spatial (and sometimes temporal) field inversion that augments closures while enforcing parsimony.
\citet{parish2016} introduced FIML's two-stage procedure: (i) a partial differential equation (PDE)-constrained inversion to recover a correction field by minimizing loss in available data plus regularization, followed by (ii) supervised learning to map flow features to the inferred correction.
\citet{singh2016} demonstrated functional-error inversion with Bayesian formulations, targeting force coefficients and surface distributions across multiple flow conditions.
For unsteady flows, \citet{fang2024} proposed time-accurate FIML where the loss aggregates spatio-temporal errors, and recent variants enforce physical invariance constraints using ensemble-Kalman regularization \citep{chen2025}.
Across these works, the objective function typically targets integrated force quantities, surface distributions, or velocity profiles---quantities that are either low-dimensional summaries or confined to specific spatial regions.

Extending field inversion to full spatio-temporal state matching, however, faces fundamental obstacles.
The objective space becomes astronomically high-dimensional, scaling with the product of spatial degrees of freedom and temporal snapshots.
Optimization landscapes grow increasingly ill-conditioned, plagued by local minima and sensitivity to initialization.
Computational costs become prohibitive as resolution increases, and point-wise field matching provides no inherent regularization against overfitting to noise or measurement artifacts.
These challenges are not merely practical inconveniences---they represent a fundamental mismatch between the structure of the inverse problem and the representation used to solve it.

Modal analysis offers a path forward~\cite{Taira2017}.
By decomposing complex spatio-temporal fields into a hierarchy of coherent structures, modal techniques provide principled dimensionality reduction that preserves essential physics while yielding well-conditioned, low-dimensional representations.
Among modal methods, POD is particularly well-suited for inverse problems: it provides an optimal low-rank basis that captures dominant flow structures, and its data-driven nature requires no prior knowledge of the underlying dynamics.

POD is a statistical technique that extracts orthogonal modes capturing the dominant energy content of complex systems~\citep{Berkooz1993}.
It is equivalent to the Karhunen-Lo\`eve decomposition and principal component analysis (PCA), providing an optimal low-rank basis in the least-squares sense \citep{Holmes1996}.
\citet{Sirovich1987} introduced the snapshot method, which made POD computationally feasible for large-scale fluid datasets.
Since then, POD has become a cornerstone of modal decomposition in fluid mechanics \citep{Rowley2005, Taira2017}.
POD enables coherent structure extraction and reduced-order modeling across multiple engineering disciplines \citep{Noack2003, Benner2015}.
In unsteady flows, POD-based Galerkin models capture transient and post-transient wake dynamics \citep{Noack2003}.
Balanced POD (BPOD) improves reduced-order-model (ROM) controllability and observability by embedding control-theoretic concepts \citep{Rowley2005}.
POD-ROMs enable efficient flow reconstruction and parametric analysis in aerodynamic design and inverse problems, with robust formulations addressing numerical stability and parameter dependence \citep{BuiThanh2004, Bergmann2009}.
Applications extend to acoustics, thermo-fluids, and reacting flows, delivering substantial computational savings in optimization and parametric studies \citep{Benner2015}.

POD has been widely used in gradient-based optimization and control.
POD-based optimization methods span Hamilton--Jacobi--Bellman (HJB)-based reduced-state control \citep{Kunisch2004, Kunisch2008}, surrogate-based inverse design with fixed bases \citep{BuiThanh2004}, adjoint-driven POD-ROM optimization \citep{Bergmann2008, Yamaleev2010}, analytic parametric sensitivity-based POD enrichment \citep{Akkari2014}, POD-assisted structural and topology optimization \citep{Ferro2019, Karcher2022}, and weighted or differentiable POD formulations with explicit basis sensitivities \citep{Merle2023, VanSchie2025}.

For instance, \citet{Kunisch2004, Kunisch2008} coupled POD with HJB feedback control, optimizing reduced states via dynamic programming, with gradients derived from the reduced HJB system rather than POD-basis derivatives.
In aerodynamic inverse design, \citet{BuiThanh2004} employed POD surrogates within gradient-based loops, computing sensitivities with respect to surrogate coefficients while keeping the POD basis fixed.
Similarly, for laminar cylinder-wake control, \citet{Bergmann2008} minimized drag using a trust-region POD-ROM with adjoint-based gradients on the reduced system, and \citet{Yamaleev2010} introduced a local-in-time discrete adjoint for unsteady optimization, treating the POD basis as frozen.
Thus, most POD-based optimizations operate on reduced states while treating the POD basis as fixed (frozen-mode optimization).

Explicit POD-basis differentiation remains rare.
\citet{Akkari2014} derived analytic solution sensitivities for parametric parabolic equations and incorporated them into POD basis enrichment, without explicitly differentiating the POD modes.
In structural and topology optimization, POD has been used to reduce design spaces, but gradients were taken with respect to classical density or topology variables on fixed or periodically updated bases \citep{Ferro2019, Karcher2022}.
\citet{Merle2023} combined POD geometric sensitivities with adjoint flow and mesh-deformation gradients, without differentiating the flow-induced POD basis, while \citet{VanSchie2025} proposed weighted POD with explicit reduced-basis derivatives embedded in a gradient-based optimizer.

A key challenge in differentiable modal analysis is the tight coupling between modal decomposition and the underlying dynamical system: the flow solution parametrizes the modal analysis itself.
This dependence arises in linear stability analysis through steady base flows \citep{Shi2019a, He2021c}, in bifurcation theory via Lyapunov-coefficient equations parametrized by steady states \citep{Kuznetsov2023, He2025a}, and in data assimilation \citep{Wu2025a}.

Derivative computation methods include finite differences (FD), complex step (CS), algorithmic differentiation (AD), and analytic direct and adjoint approaches, each trading accuracy, cost, and scalability \citep[Chapter~6]{Martins2022}.
FD and CS are black-box and easy to deploy, with CS achieving machine precision, but both scale poorly with the number of inputs, motivating adjoint-type methods for high-dimensional problems \citep{Martins2003a, Lyu2014f, Martins2013a, Martins2022}.
AD applies the chain rule through the computation graph: forward mode scales with inputs, reverse mode with outputs, but incurs significant storage costs \citep{Martins2022}, with recent reverse-mode AD formulations developed for eigenvalue problems \citep{He2020a, He2021b}.
Analytic direct and adjoint methods are preferable when inputs are fewer or greater than outputs, respectively \citep[Chapter~6]{Martins2022}. Recent work has differentiated the singular value decomposition (SVD) using adjoints \citep{Kanchi2025a}, and the present work leverages a differentiable SVD-POD connection \citep{Taira2017}.

To address these challenges, we propose MCFI, a paradigm that reformulates the inverse problem in the reduced space of POD modes rather than the full physical state space (see \Cref{sec:mcfi} and \Cref{fig:mcfi_schematic}).
The core idea is simple: \textbf{\emph{match structures, not points}}.
By targeting dominant flow structures instead of point-wise field values, MCFI provides a compact, physically meaningful objective that naturally regularizes the inversion and dramatically reduces computational burden.
Coherent structures captured by leading POD modes encode the essential dynamics; matching these modes implicitly enforces global flow consistency without requiring explicit spatio-temporal point matching.

Central to MCFI is the need for a differentiable POD framework.
We develop an adjoint-based method that efficiently computes sensitivities of POD modes with respect to model parameters, enabling gradient-based optimization in the modal space.
This framework handles the tight coupling between the dynamical system and the modal decomposition, and \textbf{requires only one unsteady adjoint solve regardless of the number of POD modes in the objective function}.
We demonstrate MCFI on one- and two-dimensional modified viscous Burger's equations, optimizing spatially varying coefficients to match target dynamics through mode-matching.

The remainder of this paper is organized as follows.
\Cref{sec:gov_eq} presents the coupled unsteady dynamical system and POD governing equations.
\Cref{sec:derivatives} develops the adjoint-based differentiable framework for multi-mode POD.
\Cref{sec:mcfi} introduces the MCFI strategy.
\Cref{sec:resul_dsc} applies the proposed methods to modified Burger's equations, first verifying the adjoint gradients against finite differences, then demonstrating adjoint-based optimization for 1D and 2D MCFI problems.
\Cref{sec:conc} presents our conclusions.

\section{Coupled dynamical system and POD analysis}\label{sec:gov_eq}

In this section, we present the general form of governing equations we considered in current study and the governing equations for POD.
In the first part, we discuss the general form of unsteady equation of the governing equations in residual-state form.
In the second subsection, we discuss the undertaken POD computation approach leveraging the relationship between SVD and POD, and present the POD governing equation in residual-state form as well.

\subsection{Dynamical system governing equation}\label{sec:gov_eq_base}
In the current study, we consider dynamical systems of the form
\begin{equation}
\f{\d \mb{u}}{\d t} = \mb{r}(\mb{u}, \mb{x}),
\end{equation}
where $\mb{r}\in\mathbb{R}^{n_s}$ is the dynamical-system spatial residual vector, $\mb{u}\in\mathbb{R}^{n_s}$ is the dynamical-system state vector, and $\mb{x}\in\mathbb{R}^{n_x}$ is the design variables vector.
Here, $n_s$ is the number of spatial degrees of freedom, or simply the number of grid points post discretization where the state variables and residuals are evaluated at, and $n_x$ is the total number of inputs or design variables.

Different time integration schemes can be adopted to solve the unsteady problem, e.g., when we apply the forward Euler integration scheme, we have
\begin{equation}\label{eq:forward_euler}
\mbg{u}^{\,(i+1)} \;=\; \mbg{u}^{\,(i)} \;+\; \Delta t\,\mb{r}(\mbg{u}^{\,(i)},\mb{x}),
\end{equation}
where $i$ is the time-step index and $\Delta t$ is the time-step size.
\Cref{eq:forward_euler} can be expanded following the general forward Euler style in residual-state form
\begin{equation}
\mb{r}_{\mathrm{uns}}(\mb{u}, \mbg{\mb{x}}) = 
\begin{bmatrix}
\mb{u}^{(1)} \\
\mb{u}^{(2)} \\
\vdots \\
\mb{u}^{(n_t)} \\
\end{bmatrix}
-
\begin{bmatrix}
\mb{u}^{(0)} + \Delta t \text{ } \mb{r}(\mb{u}^{(0)}, \mbg{\mb{x}}) \\
\mb{u}^{(1)} + \Delta t \text{ } \mb{r}(\mb{u}^{(1)}, \mbg{\mb{x}}) \\
\vdots \\
\mb{u}^{(n_t - 1)} + \Delta t \text{ } \mb{r}(\mb{u}^{(n_t - 1)}, \mbg{\mb{x}}) \\
\end{bmatrix}, 
\quad 
\mb{u} = 
\begin{bmatrix}
\mb{u}^{(1)} \\
\mb{u}^{(2)} \\
\vdots \\
\mb{u}^{(n_t)} \\
\end{bmatrix},\label{eq:rbase}
\end{equation}
where $\mb{r}_{\mathrm{uns}}\in\mathbb{R}^{n_s n_t}$ is the dynamical-system residual, $\mb{u}\in\mathbb{R}^{n_s n_t}$ is the stacked state vector, $\mb{r}$ is the spatial residual, and $n_t$ is the number of stored time levels (snapshots) in the unsteady trajectory.
This summarizes the governing equation general form for the dynamical system.

\subsection{Proper orthogonal decomposition governing equation}\label{sec:gov_eq_pod}

POD computation can be achieved through the method of snapshots which solves the eigenvalue problem (EVP), direct SVD of the snapshot matrix whose POD is sought, or solving the EVP for the covariance matrix~\cite{Taira2017}.
In the current manuscript, the direct SVD approach was undertaken.
In lieu of this, we state the governing equations to the SVD and a short introduction to the governing equation of the symmetric embedding matrix method (SEMM) proposed by Kanchi and He~\cite{Kanchi2025a}, which uses a special form of the SVD governing equation.

SVD matrix factorization for a general real-valued snapshot matrix $\tilde{\mb{U}}\in\mathbb{R}^{n_s \times n_t}$ is defined as 
\begin{equation}
\tilde{\mb{U}}=\mbg{\Phi\Sigma}\mb{V}^\intercal,
\label{eq:SVD_gen_equn}
\end{equation}
where $\mbg{\Phi}\in\mathbb{R}^{n_s \times n_s}$ collects the left singular vectors, $\mbg{\Sigma}\in\mathbb{R}^{n_s \times n_t}$ is the rectangular diagonal matrix of singular values, and $\mb{V}\in\mathbb{R}^{n_t \times n_t}$ collects the right singular vectors.
Here $\tilde{\mb{U}}$ is the zero-mean snapshot matrix obtained by removing the temporal mean from each row of the snapshot matrix $\mb{U}=[\mb{u}^{(1)}, \ldots, \mb{u}^{(n_t)}]$ (see \Cref{eq:centering_snap}).

The POD governing equation is
\begin{equation}
\begin{aligned}
\tilde{\mb{U}}\mb{v} &= \sigma \mbg{\phi},\\
\tilde{\mb{U}}^\intercal\mbg{\phi} &= \sigma \mb{v},\\
\mbg{\phi}^\intercal \mbg{\phi} &= 1,
\end{aligned} 
\label{eq:svd_gov_real}
\end{equation}
where $\{\mbg{\phi},\sigma,\mb{v},\tilde{\mb{U}}\} \in \mathbb{R}$.
We refer the readers to paper by Kanchi and He~\cite{Kanchi2025a} which presents the derivation of these governing equations for a general complex matrix $\mb{A}\in\mathbb{R}^{n_s \times n_t}$.

We can now write this in the residual and state variables form as follows
\begin{equation} 
\mb{r}_{\mathrm{POD},i} = 
\begin{bmatrix}
\tilde{\mb{U}}_i \mb{v}_i - \sigma_i \mbg{\phi}_i,\\
\tilde{\mb{U}}_i^\intercal\mbg{\phi}_i - \sigma_i \mb{v}_i,\\
\mbg{\phi}_i^\intercal \mbg{\phi}_i - 1\\
\end{bmatrix}
,\text{ }
\mb{w}_i =
\begin{bmatrix}
\mbg{\phi}_i\\
\mb{v}_i\\
\sigma_i
\end{bmatrix},
\label{eq:rmodal}
\end{equation}
where $\mb{w}_i\in\mathbb{R}^{n_s+n_t+1}$ is the modal state vector and $\mb{r}_{\mathrm{POD},i}\in\mathbb{R}^{n_s+n_t+1}$ is the modal residual for the $i$-th mode.

\subsection{Global residual form}\label{sec:global_res_form}

In the preceding subsections, we established the dynamical-system and modal residual and state vectors.
\Cref{eq:rbase} gives the dynamical-system residual $\mb{r}_{\mathrm{uns}}$ and state $\mb{u}$.
\Cref{eq:rmodal} gives the modal residual $\mb{r}_{\mathrm{POD}} $ and modal state $\mb{w}$ for a single mode.

When multiple POD modes are of interest, each mode $i\in\{1,\ldots,m\}$ introduces its own singular triplet $(\mbg{\phi}_i, \mb{v}_i, \sigma_i)$ and corresponding modal state $\mb{w}_i$ and residual $\mb{r}_{\mathrm{POD},i}$, where $m$ is the number of POD modes retained.
The general global residual and state form of the equations can be formed by stacking the dynamical system together with all $m$ modal systems:
\begin{equation}
\mb{r}_\mathrm{global}(\mb{q}; \mb{x}) =
\begin{bmatrix}
\mb{r}_{\mathrm{uns}}(\mb{u}; \mb{x}) \\
\mb{r}_{\mathrm{POD},1}(\mb{w}_1; \mb{u}) \\
\vdots \\
\mb{r}_{\mathrm{POD},m}(\mb{w}_m; \mb{u}) \\
\end{bmatrix},
\quad
\mb{q} =
\begin{bmatrix}
\mb{u} \\
\mb{w}_1 \\
\vdots \\
\mb{w}_m \\
\end{bmatrix},
\label{eq:global_multi-mode}
\end{equation}
where $\mb{q}\in\mathbb{R}^{n_s n_t + m(n_s+n_t+1)}$ is the global state vector and $\mb{r}_\mathrm{global}\in\mathbb{R}^{n_s n_t + m(n_s+n_t+1)}$ is the global residual vector.
Each modal residual $\mb{r}_{\mathrm{POD},i}$ depends on the dynamical-system state $\mb{u}$ through the snapshot matrix but is independent of the other modal states $\mb{w}_j$ for $j\neq i$.
This feedforward structure, where the dynamical system drives all modal systems without inter-modal coupling, is key to the efficient adjoint formulation presented in~\Cref{sec:derivatives}.

As a special case, when only a single mode ($m=1$) is of interest, the global system reduces to
\begin{equation}
\mb{r}_\mathrm{global}(\mb{q}; \mb{x}) =
\begin{bmatrix}
\mb{r}_{\mathrm{uns}}(\mb{u}; \mb{x}) \\
\mb{r}_{\mathrm{POD}}(\mb{w}; \mb{u}) \\
\end{bmatrix},
\quad
\mb{q} =
\begin{bmatrix}
\mb{u} \\
\mb{w} \\
\end{bmatrix}.
\label{eq:global_singlemode}
\end{equation}
This simplified form is used in the following derivations for clarity; the extension to multiple modes is straightforward.

\section{Function of interest}
For a general function of interest that may depend on the dynamical-system flow and one or more POD modes,
\begin{equation}
f = f(\mb{q}; \mb{x}) = f(\mb{u}, \mb{w}_1, \ldots, \mb{w}_m; \mb{x}).
\end{equation}
The framework accommodates objectives involving the dynamical-system state $\mb{u}$, the modal states $\mb{w}_i = (\mbg{\phi}_i, \mb{v}_i, \sigma_i)$, or both.
In this work, we demonstrate mode-matching objectives of the form
\begin{equation}
\begin{aligned}
f &= \tfrac{1}{2}\big\|\mbg{\phi} - \mbg{\phi}^\star\big\|_2^2, &
&\text{(quadratic mode loss)}\\
f &= \big\|\mbg{\phi} - \mbg{\phi}^\star\big\|_2 + (\sigma - \sigma^\star)^2, &
&\text{(mode and energy loss)}\\
f &= \sum_{i=1}^{m} \Big(\big\|\mbg{\phi}_i - \mbg{\phi}_i^\star\big\|_2 + (\sigma_i - \sigma_i^\star)^2\Big), &
&\text{(multiple modes and energies loss)}
\end{aligned}
\label{eq:foi_examples}
\end{equation}
where $\mbg{\phi}_i^\star$ and $\sigma_i^\star$ denote target mode shapes and singular values, respectively.
Beyond mode matching, the differentiable POD framework naturally extends to objectives that couple the dynamical-system and modal quantities.
Classic field inversion (FI) targets mean-flow matching:
\begin{equation}
f = \big\|\bar{\mb{u}} - \bar{\mb{u}}^\star\big\|_2^2, \quad \text{(FI: mean flow matching)}
\label{eq:foi_fi}
\end{equation}
where $\bar{\mb{u}}$ is the temporal mean of the dynamical-system flow.
MCFI introduces new objectives that incorporate modal information:
\begin{equation}
\begin{aligned}
f &= \big\|\bar{\mb{u}} - \bar{\mb{u}}^\star\big\|_2^2 + \lambda \sum_{i=1}^{m} \big\|\mbg{\phi}_i - \mbg{\phi}_i^\star\big\|_2^2, && \text{(mean flow + mode matching)}\\
f &= g(\mb{u}) + \lambda \sum_{i=1}^{m} \sigma_i, && \text{(dynamical objective + energy penalty)}\\
f &= -\f{\sigma_1}{\sigma_2}, && \text{(spectral gap maximization)}
\end{aligned}
\label{eq:foi_extensions}
\end{equation}
where $g(\mb{u})$ is any dynamical-system objective (e.g., drag, lift) and $\lambda$ is a weighting parameter.
The adjoint-based derivative computation strategy presented in \Cref{sec:derivatives} applies to any objective of the form $f(\mb{u}, \mb{w}_1, \ldots, \mb{w}_m; \mb{x})$, including all examples above. 

\section{Differentiable coupled dynamical system and POD analysis}\label{sec:derivatives}
In this section, we present the derivative computation methodology for the differentiable POD using the adjoint method.
In the first subsection, we present the overview of the derivative computation workflow by splitting it into two parts, adjoint based differentiable SVD and the unsteady adjoint for the dynamical system derivative, with the motivation for doing it.
Next, we present the differentiable SVD methodology and finally present the unsteady adjoint methodology used for the differentiable modal analysis of a general dynamical system.

\subsection{Differentiable modal analysis}
We apply the adjoint method to compute the derivative
\begin{equation}
\f{\d f}{\d \mb{x}} = \f{\p f}{\p \mb{x}} - \boldsymbol{\psi}^\intercal \f{\p \mb{r}}{\p \mb{x}}, \quad \f{\p \mb{r}}{\p \mb{q}}^\intercal\boldsymbol{\psi} = \f{\p f}{\p \mb{q}}^\intercal,
\label{eq:Adjoint equation}
\end{equation}
where $\boldsymbol{\psi}$ is the adjoint variables vector and the second equation is the adjoint equation~\cite{Martins2022}.
This adjoint system encompasses the residual and state from the base dynamical system~\Cref{eq:rbase} and its modal system~\Cref{eq:rmodal}.

When multiple POD modes contribute to the objective, the same pipeline applies to each singular triplet before summing the modal forcings.
\begin{equation}
\begin{bmatrix}
\f{\p \mb{r}_{\mathrm{uns}}}{\p \mb{u}}^\intercal & \f{\p \mb{r}_{\mathrm{POD}, 1} }{\p \mb{u}}^\intercal  & \ldots & \f{\p \mb{r}_{\mathrm{POD}, m} }{\p \mb{u}}^\intercal \\
\bcancel{\f{\p \mb{r}_{\mathrm{uns}}}{\p \mb{w}_1}^\intercal}{(=\mb{0})} & \f{\p \mb{r}_{\mathrm{POD}, 1} }{\p \mb{w}_1}^\intercal & \ldots & \bcancel{\f{\p \mb{r}_{\mathrm{POD}, 1}}{\p \mb{w}_m}^\intercal}{(=\mb{0})} \\
\vdots & \vdots & \ddots & \vdots \\
\bcancel{\f{\p \mb{r}_{\mathrm{uns}}}{\p \mb{w}_m}^\intercal}{(=\mb{0})} & \bcancel{\f{\p \mb{r}_{\mathrm{POD}, 1}}{\p \mb{w}_m}^\intercal}{(=\mb{0})}  & \ldots & \f{\p \mb{r}_{\mathrm{POD}, m} }{\p \mb{w}_m}^\intercal \\
\end{bmatrix}
\begin{bmatrix}
\boldsymbol{\psi}_{\mathrm{uns}} \\
\boldsymbol{\psi}_{\mathrm{POD}, 1} \\
\vdots \\
\boldsymbol{\psi}_{\mathrm{POD}, m}
\end{bmatrix}
= 
\begin{bmatrix}
\f{\p f}{\p \mb{u}}^\intercal \\
\f{\p f}{\p \mb{w}_1}^\intercal \\
\vdots \\
\f{\p f}{\p \mb{w}_m}^\intercal \\
\end{bmatrix}.
\label{eq:adj_block_sparse}
\end{equation}
Here, the sparsity pattern of the global residual Jacobian $\p \mb{r}_\mathrm{global} / \p \mb{q}$ is evident (see the matrix on the left-hand side of~\Cref{eq:adj_block_sparse}). 
Leveraging this sparsity pattern, we can very efficiently compute the derivative without directly solving the large coupled adjoint equation. 
\textbf{No matter how many modes of interest, we only need to solve one unsteady adjoint equation}.



For clarity, we present the single-mode case ($m=1$); the multi-mode extension follows the same block structure with additional modal rows.
The adjoint equation can be expanded as
\begin{equation}
\begin{bmatrix}
\f{\p \mb{r}_{\mathrm{uns}}}{\p \mb{u}}^\intercal & \f{\p \mb{r}_{\mathrm{POD}} }{\p \mb{u}}^\intercal \\
\bcancel{\f{\p \mb{r}_{\mathrm{uns}}}{\p \mb{w}}^\intercal}{(=\mb{0})} & \f{\p \mb{r}_{\mathrm{POD}} }{\p \mb{w}}^\intercal \\
\end{bmatrix}
\begin{bmatrix}
\boldsymbol{\psi}_{\mathrm{uns}} \\
\boldsymbol{\psi}_{\mathrm{POD}}
\end{bmatrix}
= 
\begin{bmatrix}
\f{\p f}{\p \mb{u}}^\intercal \\
\f{\p f}{\p \mb{w}}^\intercal
\end{bmatrix}.
\label{eq:adj_block}
\end{equation}
The off-diagonal zero block arises because the unsteady residual does not depend on the POD variables $\mb{w}$.
This feedforward structure allows the system to be solved via back-substitution:
\begin{equation}
\label{eq:adj}
\begin{aligned}
\boldsymbol{\psi}_{\mathrm{POD}} &\leftarrow  \f{\p \mb{r}_{\mathrm{POD}} }{\p \mb{w}}^\intercal \boldsymbol{\psi}_{\mathrm{POD}} = \f{\p f}{\p \mb{w}}^\intercal,\\
\boldsymbol{\psi}_{\mathrm{uns}} &\leftarrow  \f{\p \mb{r}_{\mathrm{uns}}}{\p \mb{u}}^\intercal \boldsymbol{\psi}_{\mathrm{uns}} = \f{\p f}{\p \mb{u}}^\intercal - \f{\p \mb{r}_{\mathrm{POD}} }{\p \mb{u}}^\intercal \boldsymbol{\psi}_{\mathrm{POD}}. \\
\end{aligned}
\end{equation}
The resulting two blocks of the adjoint equation are solved separately and backsubstituted to obtain the adjoint variables.
This reduction of the large adjoint system into smaller coupled adjoint systems is more efficient.

\subsection{Differentiable proper orthogonal decomposition}\label{sec:diff_pod}
Consider the residual-state form of the POD governing equation in~\Cref{eq:rmodal} and the multi-mode formulation in~\Cref{eq:global_multi-mode}.
From the first equation in~\Cref{eq:adj}, for a single mode $i$
\begin{equation}
\f{\p \mb{r}_{\mathrm{POD},i}}{\p \mb{w}_i} = 
\begin{bmatrix}
-\sigma_i\,\mb{I}_{n_s} & \tilde{\mb{U}} & -\mbg{\phi}_i \\
\tilde{\mb{U}}^\intercal & -\sigma_i\,\mb{I}_{n_t} & -\mb{v}_i \\
2\,\mbg{\phi}_i^\intercal & \mb{0}^\intercal & 0
\end{bmatrix}, \quad \f{\p f}{\p \mb{w}_i} = \begin{bmatrix}
\f{\p f}{\p \mbg{\phi}_i}\\
\f{\p f}{\p \mb{v}_i}\\
\f{\p f}{\p \sigma_i}\\
\end{bmatrix}, \quad 
\mbg{\psi}_{\mathrm{POD},i} =
\begin{bmatrix}
\mbg{\psi}_{\mbg{\phi},i}\\
\mbg{\psi}_{\mb{v},i}\\
\psi_{\sigma,i}\\
\end{bmatrix},
\end{equation}
where $\mbg{\psi}_{\mbg{\phi},i}\in\mathbb{R}^{n_s}$, $\mbg{\psi}_{\mb{v},i}\in\mathbb{R}^{n_t}$, and $\psi_{\sigma,i}\in\mathbb{R}$ are the adjoint variables associated with the $i^\text{th}$ spatial mode, temporal coefficients, and singular-value constraint, respectively. The stacking over $i$ recovers the multi-mode form in~\Cref{eq:global_multi-mode}. 
Combining it in relation with the Jacobian of modal residual with respect to the modal state in~\Cref{eq:rmodal}, the block adjoint system for each mode $i$ is
\begin{equation}
\begin{bmatrix}
-\sigma_i\,\mb{I}_{n_s} & \tilde{\mb{U}} & -\mbg{\phi}_i \\
\tilde{\mb{U}}^\intercal & -\sigma_i\,\mb{I}_{n_t} & -\mb{v}_i \\
2\,\mbg{\phi}_i^\intercal & \mb{0}^\intercal & 0
\end{bmatrix}^\intercal 
\begin{bmatrix}
\mbg{\psi}_{\mbg{\phi},i}\\
\mbg{\psi}_{\mb{v},i}\\
\psi_{\sigma,i}\\
\end{bmatrix}
=
\begin{bmatrix}
\f{\p f}{\p \mbg{\phi}_i}\\
\f{\p f}{\p \mb{v}_i}\\
\f{\p f}{\p \sigma_i}\\
\end{bmatrix}^\intercal,
\label{eq:modal_adj_expand}
\end{equation}
where $\mb{I}_{n_s}\in\mathbb{R}^{n_s\times n_s}$ and $\mb{I}_{n_t}\in\mathbb{R}^{n_t\times n_t}$ are identity matrices matching the dimensions of $\mbg{\phi}_i$ and $\mb{v}_i$; the coefficient matrix is $(\p \mb{r}_\mathrm{POD,i}/\p \mb{w}_i)^\intercal$, the adjoint vector is $\mbg{\psi}_\mathrm{POD,i}$, and the right-hand side is $(\p f / \p \mb{w}_i)^\intercal$ in~\Cref{eq:adj}. Stacking these $m$ modal systems yields the multi-mode formulation in~\Cref{eq:global_multi-mode}.

\subsection{Unsteady adjoint for the dynamical system}\label{sec:base_adj}
Bringing our attention to~\Cref{eq:rbase,eq:adj}, we first lay out the Jacobian block, adjoint vector, and right-hand side for the dynamical-system adjoint,
\begin{equation}
\f{\p \mb{r}_{\mathrm{uns}}}{\p \mb{u}} =
\begin{bmatrix}
\mb{I} & -\mb{I} - {\Delta t}\f{\p \mb{r}_\mathrm{s}}{\p \mb{u}}\vert_{\mb{u}^{(1)}}^\intercal  &  \cdots & \mb{0} \\
\mb{0} & \mb{I} & \cdots & \mb{0}  \\
\vdots & \vdots & \ddots & \vdots\\
\mb{0} & \cdots & \cdots & \mb{I}\\
\end{bmatrix},\quad
\boldsymbol{\psi}_{\mathrm{uns}} =
\begin{bmatrix}
\mbg{\psi}_{b,1} \\
\mbg{\psi}_{b,2} \\
\vdots \\
\mbg{\psi}_{b,n_t} \\
\end{bmatrix},
\end{equation}
\begin{equation}
\f{\p f}{\p \mb{u}} =
\begin{bmatrix}
\f{\p f}{\p \mb{u}^{(1)}}\\
\f{\p f}{\p \mb{u}^{(2)}}\\
\vdots \\
\f{\p f}{\p \mb{u}^{(n_t)}}\\
\end{bmatrix}.
\end{equation}
With these definitions, the dynamical-system adjoint equation becomes
\begin{equation}
\label{eq:base_adj_expand}
\begin{bmatrix}
\mb{I} & -\mb{I} - {\Delta t}\f{\p \mb{r}_\mathrm{s}}{\p \mb{u}}\vert_{\mb{u}^{(1)}}^\intercal  &  \cdots & \mb{0} \\
\mb{0} & \mb{I} & \cdots & \mb{0}  \\
\vdots & \vdots & \ddots & \vdots\\
\mb{0} & \cdots & \cdots & \mb{I}\\
\end{bmatrix}
\begin{bmatrix}
\mbg{\psi}_{b,1} \\
\text{}\\
\mbg{\psi}_{b,2} \\
\vdots \\
\mbg{\psi}_{b,n_t} \\
\end{bmatrix}
=
\begin{bmatrix}
\f{\p f}{\p \mb{u}^{(1)}}\\
\text{}\\
\f{\p f}{\p \mb{u}^{(2)}}\\
\vdots \\
\f{\p f}{\p \mb{u}^{(n_t)}}\\
\end{bmatrix}
- \f{\p \mb{r}_{\mathrm{POD}} }{\p \mb{u}}^\intercal \boldsymbol{\psi}_{\mathrm{POD}}.
\end{equation}
We elaborate on the computation of $({\p \mb{r}_{\mathrm{POD}} }/{\p \mb{u}})^\intercal \boldsymbol{\psi}_{\mathrm{POD}}$ in~\ref{sec:modal_coupling}.
We call this term the modal-forcing term.
This step is slightly more complicated due to the presence of multiple Jacobians and the centering of the snapshot matrix in our POD workflow.

\subsection{Algorithm}

The complete procedure for computing the gradient of a modal objective with respect to design variables is summarized in~\Cref{alg:diff_pod}.
The algorithm exploits the block-diagonal structure of the modal adjoint systems, which allows each mode to be processed independently before a single unsteady adjoint solve accumulates all modal contributions.

\begin{algorithm}[H]
\caption{Differentiable POD gradient evaluation for modal objectives}
\label{alg:diff_pod}
\begin{algorithmic}[1]
\Require Snapshots $\mb{U}$, POD modes $(\mbg{\phi}_i,\sigma_i,\mb{v}_i)$ for $i=1,\ldots,m$, objective $f(\mb{u},\{\mb{w}_i\}_{i=1}^m;\mb{x})$
\Ensure Gradient $\d f/\d\mb{x}$
\State $\mb{g}\leftarrow (\partial f/\partial \mb{u})^\intercal$ \Comment{Initialize adjoint RHS}
\For{$i=1$ to $m$} \Comment{Modal adjoints}
\State $\boldsymbol{\psi}_{\mathrm{POD},i} \leftarrow \left(\partial\mb{r}_{\mathrm{POD},i}/\partial \mb{w}_i\right)^\intercal \boldsymbol{\psi}_{\mathrm{POD},i} = (\partial f/\partial \mb{w}_i)^\intercal$ \Comment{Solve, \Cref{eq:modal_adj_expand}}
\State $\mb{g} \leftarrow \mb{g} - \left(\partial\mb{r}_{\mathrm{POD},i}/\partial \mb{u}\right)^\intercal \boldsymbol{\psi}_{\mathrm{POD},i}$ \Comment{Accumulate forcing}
\EndFor
\State $\boldsymbol{\psi}_{\mathrm{uns}} \leftarrow \left(\partial\mb{r}_{\mathrm{uns}}/\partial \mb{u}\right)^\intercal \boldsymbol{\psi}_{\mathrm{uns}} = \mb{g}$ \Comment{Solve, \Cref{eq:base_adj_expand}}
\State $\d f/\d \mb{x} \leftarrow \partial f/\partial \mb{x} - \boldsymbol{\psi}_{\mathrm{uns}}^\intercal (\partial \mb{r}_{\mathrm{uns}}/\partial \mb{x})$ \Comment{Gradient, \ref{sec:residual_jacobian_wrt_x}}
\end{algorithmic}
\end{algorithm}

The key computational insight is that regardless of the number of modes $m$, \textbf{only one unsteady adjoint solve (line~6) is required}.
Each modal adjoint (line~3) is a small system involving only the singular triplet $(\mbg{\phi}_i, \mb{v}_i, \sigma_i)$, which can be solved efficiently.
The modal forcing contributions are accumulated into $\mb{g}$ before the expensive backward time integration.

\section{MCFI}\label{sec:mcfi}

\begin{figure}[H]
\centering
\includegraphics[width=1\textwidth]{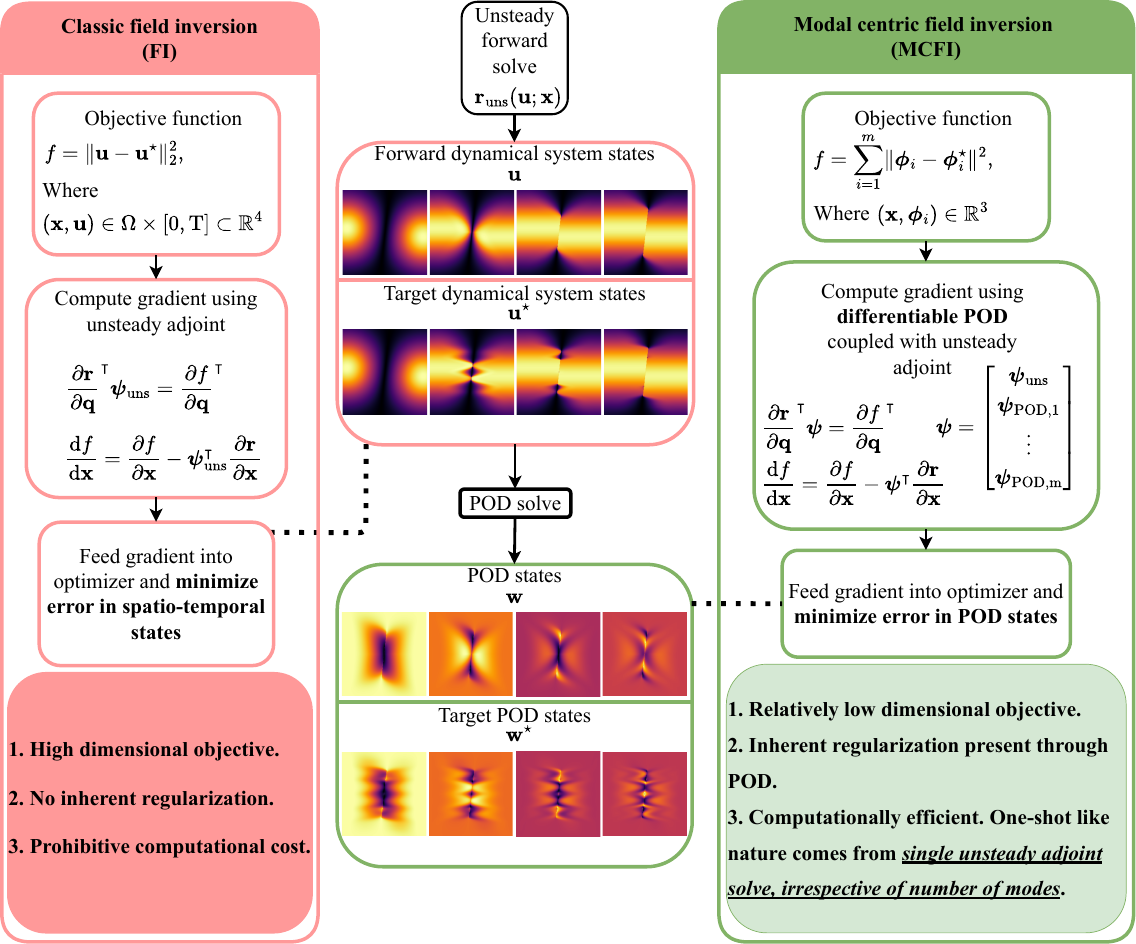}
\caption{Schematic of MCFI.
FI (left) matches point-wise values in the full spatio-temporal space, leading to high-dimensional and sometimes ill-conditioned optimization.
MCFI (right) reformulates the problem in the reduced POD state space, matching dominant coherent structures instead of individual field values.
The differentiable POD (\Cref{eq:base_adj_expand}) provides the adjoint-based gradients that enable efficient optimization in the modal space.}
\label{fig:mcfi_schematic}
\end{figure}

The differentiable POD framework enables a new class of inverse problems where the target data is expressed in terms of POD modes rather than full-field snapshots.
MCFI is illustrated in \Cref{fig:mcfi_schematic}.
Traditional field inversion seeks a spatially distributed correction field $\mb{x}$ (e.g., turbulence model discrepancy) by minimizing the loss between simulated and observed flow fields.
In contrast, MCFI operates in the reduced modal space, which offers several advantages: (i) reduced data storage since only a few dominant modes need to be matched, (ii) robustness to measurement noise concentrated in high-frequency modes, and (iii) direct physical interpretability of the target coherent structures.

\subsection{Optimization problem formulation}

The MCFI optimization problem seeks the design variables $\mb{x}$ that minimize a modal objective subject to the governing equations.
The general constrained optimization problem reads
\begin{equation}
\begin{aligned}
\min_{\mb{x}} \quad & f(\mb{u}, \mb{w}_1, \ldots, \mb{w}_m; \mb{x}) \\
\text{subject to} \quad & \mb{r}_{\mathrm{uns}}(\mb{u}; \mb{x}) = \mb{0}, \\
& \mb{r}_{\mathrm{POD},i}(\mb{w}_i; \mb{u}) = \mb{0}, \quad i = 1, \ldots, m, \\
& \mb{x}_L \leq \mb{x} \leq \mb{x}_U,
\end{aligned}
\label{eq:mcfi_opt}
\end{equation}
where $\mb{x}_L$ and $\mb{x}_U$ are lower and upper bounds on the design variables.
The equality constraints enforce the unsteady flow equations and the POD eigenvalue problems for each mode of interest.

For the mode-matching problems demonstrated in this work, the objective function takes the forms presented in~\Cref{eq:foi_examples}, where $\mbg{\phi}_i^\star$ and $\sigma_i^\star$ are the target mode shapes and singular values obtained from reference data (e.g., experimental measurements or high-fidelity simulations).
The gradient of $f$ with respect to the design variables $\mb{x}$ is computed using the adjoint method developed in~\Cref{sec:derivatives}; the total derivative formula is given in~\Cref{eq:Adjoint equation}.
The key computational advantage is that only a single unsteady adjoint solve is required regardless of the number of modes $m$ or the number of design variables in $\mb{x}$, making MCFI tractable for high-dimensional inverse problems.

\section{Results and discussion}\label{sec:resul_dsc}
In this section, we present the results for the differentiable POD framework.
We begin our analysis by introducing the governing equation for the modified two-dimensional unsteady Burgers equation.
This is followed by the forward-simulation results for the one-dimensional Burgers equation.
Next, we validate the adjoint-based derivatives by comparing them against FD computations for a selected objective in the one-dimensional setting.
We then present the optimization results for mode matching using the one-dimensional Burgers equation.
Finally, we repeat the mentioned analysis for the two-dimensional Burgers equation.


\subsection{Governing equation for the modified unsteady Burgers equation}\label{sec:burger_gov}
In this section, we describe the governing equations for the 2D modified viscous Burgers equation and present it in residual-state form.
The governing equation numerically solved in vector form is
\begin{equation}\label{eq:burgers_vector}
\partial_t\,\boldsymbol{u} \;+\; \alpha(x,y)\,\big(\boldsymbol{u}\cdot\nabla\big)\boldsymbol{u}
\;=\; \nu\,\nabla^{2}\boldsymbol{u},
\end{equation}
where $\boldsymbol{u}=[u\;\;v]^{\intercal}$, $\nu$ is viscosity and $\alpha$ is a multiplier which is the modification that we made.
The component form of~\Cref{eq:burgers_vector} is
\begin{align}
\partial_t u \;+\; \alpha(x,y)\,\big(u\,\partial_x u \;+\; v\,\partial_y u\big) &= \nu\,\big(\partial_{xx}u \;+\; \partial_{yy}u\big), \label{eq:burgers_u}\\
\partial_t v \;+\; \alpha(x,y)\,\big(u\,\partial_x v \;+\; v\,\partial_y v\big) &= \nu\,\big(\partial_{xx}v \;+\; \partial_{yy}v\big).\label{eq:burgers_v}
\end{align}
The one-dimensional Burgers equation is recovered by simply suppressing the transverse velocity and its gradients
\begin{equation}\label{eq:burgers_1d}
\partial_t u(x,t) + \alpha(x)\,u\,\partial_x u = \nu\,\partial_{xx}u,
\end{equation}
We now state the discrete spatial operators used in the solver. 
Instead of the non-differentiable $\max/\min$ switches, the code blends the backward and forward differences with a smooth relaxation weight~\cite{Harten1983}. 
For any component $w\in\{u,v\}$ at an interior node $(i,j)$ we define the smooth selectors
\begin{equation}
\omega_x(u_{i,j}) = \tfrac{1}{2}\bigl(1+\tanh(\beta\,u_{i,j})\bigr),\qquad
\omega_y(v_{i,j}) = \tfrac{1}{2}\bigl(1+\tanh(\beta\,v_{i,j})\bigr),
\end{equation}
with the steepness parameter $\beta=20$.
Further, we apply an upwind scheme for the advection term and central difference for the diffusion term. 

The spatial residual at a node $(i,j)$ in computational domain can be written compactly by collecting the nodal state as $\mbg{u}^{\,(n)}(i,j)=[\,u^{(n)}_{i,j}\; v^{(n)}_{i,j}\,]^{\intercal}$:
\begin{equation}\label{eq:spatial_residual_node}
\mbg{r}_{s}(\mbg{u}^{\,(n)})\big\vert_{i,j}
\;=\;
-\,\alpha_{i,j}\begin{bmatrix}
\big(u\,\p_x u + v\,\p_y u\big)\vert_{i,j}\\[2pt]
\big(u\,\p_x v + v\,\p_y v\big)\vert_{i,j}
\end{bmatrix}
\;+\;
\nu\begin{bmatrix}
\big(\p_{xx}u+\p_{yy}u\big)\vert_{i,j}\\[2pt]
\big(\p_{xx}v+\p_{yy}v\big)\vert_{i,j}
\end{bmatrix},
\end{equation}
where each derivative on the right is replaced by the corresponding discrete expression above.
Stacking all grid values into the global state $\mbg{u}^{\,(n)}\in\mathbb{R}^{2N}$ and defining the stacked spatial residual $\mbg{r}_{s}(\mbg{u}^{\,(n)})\in\mathbb{R}^{2N}$, the explicit time update (forward Euler) is
\begin{equation}\label{eq:update_compact}
\mbg{u}^{\,(n+1)} \;=\; \mbg{u}^{\,(n)} \;+\; \Delta t\,\mbg{r}_{s}(\mbg{u}^{\,(n)},\mb{x}),
\end{equation}
where $N=n_x \times n_y$.
This summarizes the solver's computation methodology in compact form.
The equation is also in the same form as~\Cref{eq:forward_euler}.

\subsection{1D Burgers equation}\label{sec:1D_Burgers}
In this section, we present the results for the forward run for 1D Burgers equation.
The 1D Burgers equation uses the same form of governing equation described in the preceding subsection.
We simply ignore the second ($y$) dimension.
Next, we present the verification of the adjoint gradients against those from FD and finally present results for mode-matching optimization.
\subsubsection{Results for the modified unsteady Burgers 1D equation}\label{sec:results1d}
The 1D Burgers PDE was solved numerically with the numerical scheme described in~\Cref{sec:burger_gov}.
The PDE solved is as shown in~\Cref{eq:burgers_1d}.
We use the forward Euler time integration.
This is evident from~\Cref{eq:update_compact}.
The computational domain is 1D with $x$ $\in[-1,1]$.
The initial condition superposes two Gaussian pulses with a sinusoidal background,
\begin{equation}
u(x,0) = e^{-(x+0.7)^2/0.05}
- 0.8\,e^{-(x-0.7)^2/0.05}
+ 0.25 \sin(2\pi x),
\label{eq:burgers1d_ic}
\end{equation}
where $x \in [-1,1]$.
The spatial grid uses $161$ uniform points on $x\in[-1,1]$ with viscosity $\nu=5\times 10^{-4}$.
The smooth upwind blending (Section~\ref{sec:burger_gov}) allows us to march explicitly with a fixed time step determined by the Courant-Friedrichs-Lewy (CFL) number's restriction $\Delta t = 0.3\,\Delta x / \max |u|$ based on the initial condition.
For the present case this yields $\Delta t\approx5.0\times 10^{-3}$, and the simulation was advanced to $t_{\mathrm{end}}=2.5$ s using $496$ steps.
No adaptive re-scaling of the CFL was required because the smooth fluxes keep the Jacobian well behaved even as the wave steepens.

\begin{figure}[H]
\centering
\includegraphics[width=0.85\linewidth]{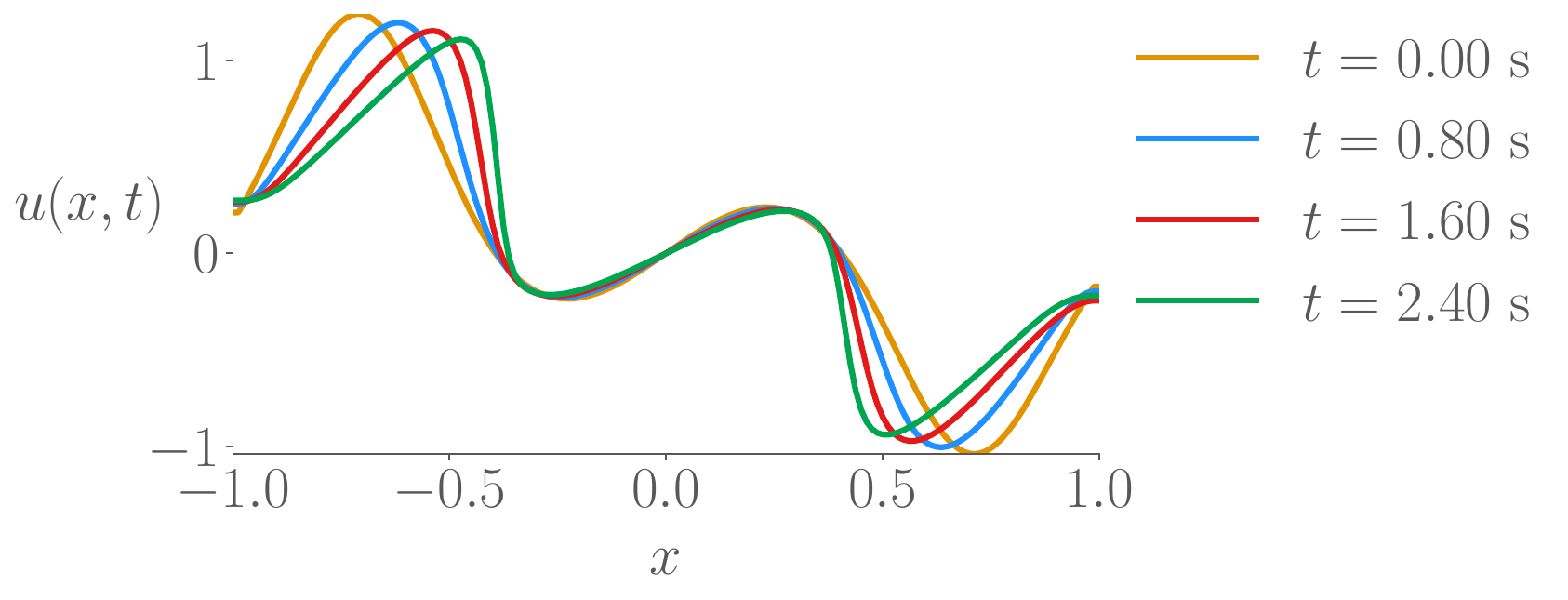}
\caption{Selected snapshots of the 1D viscous Burgers solution for the uniform multiplier case ($\alpha(x)=0.1$), with the legend denoting $t=\{0.00,0.80,1.60,2.40\}$ . 
The smooth upwind/central discretization captures the pulse interaction and subsequent viscous decay resolved in the POD snapshots.}
\label{fig:burgers1d_snapshots}
\end{figure}

\Cref{fig:burgers1d_snapshots} shows four snapshots of the baseline simulation driven by the uniform multiplier $\alpha(x)=0.1$, with physical run time of $2.50$ s to capture the full transient.
The initial pair of Gaussian pulses steepens, interacts with the boundary forcing, and ultimately diffuses away; the legend tracks the profiles at $t=\{0.00,\,0.80,\,1.60,\,2.40\}$ s, which span the energetic portion of the run that excites the leading POD mode used for mode matching.
The progressive flattening of the leading crest around $x\approx-0.5$ and the sharpening of the trailing gradient near $x\approx0.5$ illustrate the expected Burgers shock-like behavior that the differentiable solver captures without spurious oscillations.

\subsubsection{Comparison of FD and adjoint method}\label{sec:adj1d}
We applied the adjoint method presented in~\Cref{sec:derivatives}.
The finite-difference gradients used for verification rely on the forward-difference stencil summarized in~\ref{sec:fd_verification}.
The design variables are four Gaussian control coefficients $\boldsymbol{\alpha}_c=[\alpha_{c,1},\allowbreak\alpha_{c,2},\allowbreak\alpha_{c,3},\allowbreak\alpha_{c,4}]^\intercal$ that scale localized bumps centered at $x=\{-0.70,\,-0.15,\,0.40,\,0.75\}$ with width $\Delta x_c=0.25$.
Their contribution to the advection coefficient is
\begin{equation}
\alpha(x;\boldsymbol{\alpha}_c) = 0.1 + \sum_{j=1}^{4} \alpha_{c,j}\,b_j(x), \qquad -0.35 \le \alpha_{c,j} \le 0.35,
\label{eq:alpha_ctrl_1d}
\end{equation}
where each $b_j(x)$ is a Gaussian basis function centered at $x_{c,j}$,
\begin{equation}
b_j(x) = e^{-\left(\frac{x-x_{c,j}}{\Delta x_c}\right)^2},
\end{equation}
with width parameter $\Delta x_c = 0.25$.
The basis functions are evaluated on the discrete solver grid $\{x_i\}_{i=1}^{161}$.
The baseline simulation starts from $\boldsymbol{\alpha}_c^{(0)}=\mb{0}$ so that $\alpha(x;\boldsymbol{\alpha}_c^{(0)})=0.1$, while the target multiplier is formed from $\boldsymbol{\alpha}_c^\star=[0.25,\,-0.15,\,0.05,\,0.15]^\intercal$, yielding $\alpha^\star(x)=\alpha(x;\boldsymbol{\alpha}_c^\star)$.

The resulting multiplier reads as in~\Cref{eq:alpha_ctrl_1d}, and the objective penalizes the loss between the leading POD mode and a prescribed target extracted from the reference control:
\begin{equation}\label{eq:burgers1d_objective}
f=\tfrac{1}{2}\big\|\mbg{\phi}_1(\boldsymbol{\alpha})-\mbg{\phi}_1^\star\big\|_2^2,
\end{equation}
where $\mbg{\phi}_1$ is the first left singular vector of the zero-mean snapshot matrix and $\mbg{\phi}_1^\star$ is is target leading POD mode generated by one forward run with control $\boldsymbol{\alpha} = \boldsymbol{\alpha}_c^\star$.
The adjoint and FD sensitivity evaluations follow exactly the coupled formulation summarized in Sec.~\ref{sec:derivatives} (see especially~\Cref{eq:Adjoint equation,eq:modal_adj_expand,eq:base_adj_expand}).

Table~\ref{tab:1d_grad_check} reports the gradient verification for the test point $\boldsymbol{\alpha}_c=[0.90,\allowbreak -0.15,\allowbreak 0.05,\allowbreak 0.15]^\intercal$.
The adjoint and FD gradients agree to within $O(10^{-6})$ in absolute value and tighter than $2.5\times 10^{-5}$ in relative error across all components, confirming that the 1D implementation faithfully reproduces the general adjoint derivatives.

\begin{table}[H]
\centering
\ra{1.2}
\caption{Adjoint versus finite-difference gradients for the 1D Burgers Gaussian controls used in the verification case.}
\label{tab:1d_grad_check}
\vspace{1ex}
\resizebox{\columnwidth}{!}{%
\begin{tabular}{l r r r r}
\toprule
\multicolumn{1}{c}{Component} &
\multicolumn{1}{c}{Adjoint gradient} &
\multicolumn{1}{c}{FD gradient} &
\multicolumn{1}{c}{$\vert\Delta\vert$} &
\multicolumn{1}{c}{Relative error} \\
\midrule
$\alpha_{c,1}$ & $\phantom{-}0.2226819780$ & $\phantom{-}0.222681\underline{8384}$ & $1.40\times10^{-7}$ & $6.29\times10^{-7}$ \\
$\alpha_{c,2}$ & $\phantom{-}0.1103718593$ & $\phantom{-}0.110371\underline{8646}$ & $5.33\times10^{-9}$ & $4.83\times10^{-8}$ \\
$\alpha_{c,3}$ & $-0.1575599861$ & $-0.15755\underline{62218}$ & $3.76\times10^{-6}$ & $2.39\times10^{-5}$ \\
$\alpha_{c,4}$ & $-0.1433312851$ & $-0.1433\underline{295612}$ & $1.72\times10^{-6}$ & $1.21\times10^{-5}$ \\
\bottomrule
\end{tabular}%
}
\end{table}

\subsubsection{Single mode MCFI for 1D Burger's equation and optimization}\label{sec:1_mode_mcfi}
\begin{figure}[H]
\centering
\includegraphics[width=0.95\linewidth]{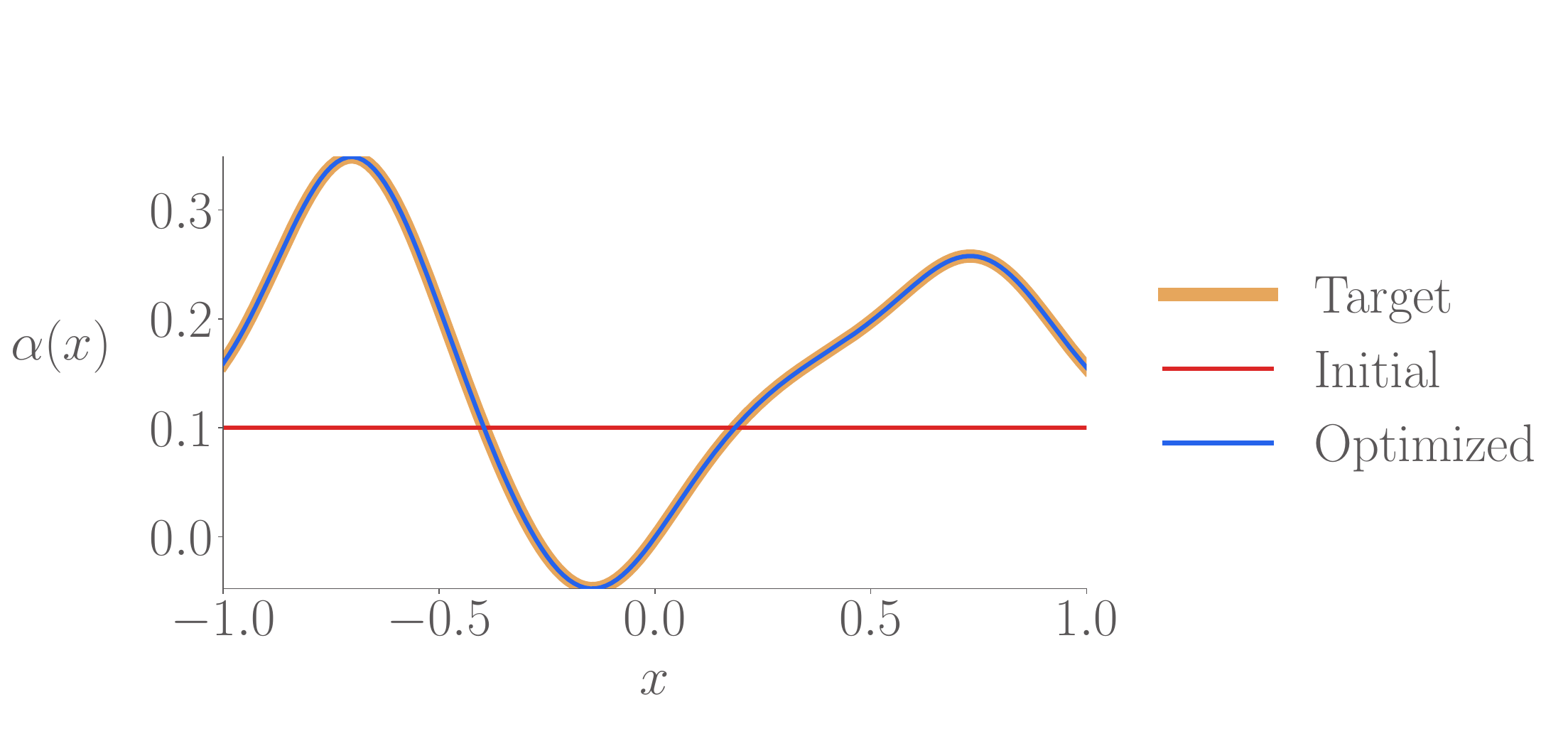}
\caption{Comparison of the initial, optimized, and target Gaussian control profiles for the 1D mode-matching problem. The optimized profile recovers the target everywhere while preserving the smoothing enforced by the Gaussian basis.}
\label{fig:burgers1d_alpha_cmp}
\end{figure}

\begin{figure}[H]
\centering
\includegraphics[width=0.95\linewidth]{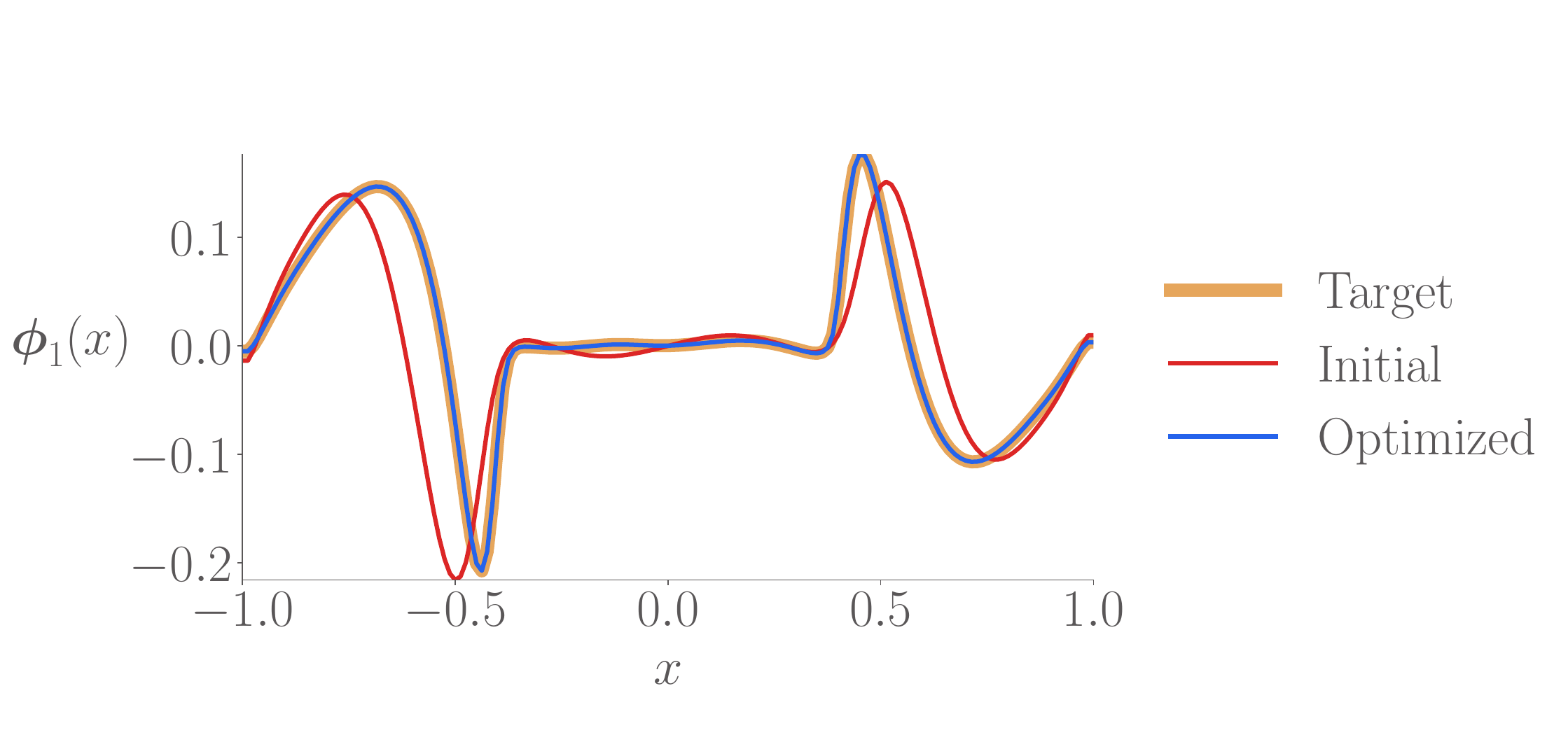}
\caption{First POD mode associated with the initial control, the optimized solution, and the target trajectory. The optimized mode collapses onto the target mode, verifying that the POD-adjoint coupling transports the sensitivities without loss of accuracy.}
\label{fig:burgers1d_mode_cmp}
\end{figure}

\begin{figure}[H]
\centering
\includegraphics[width=1.0\linewidth]{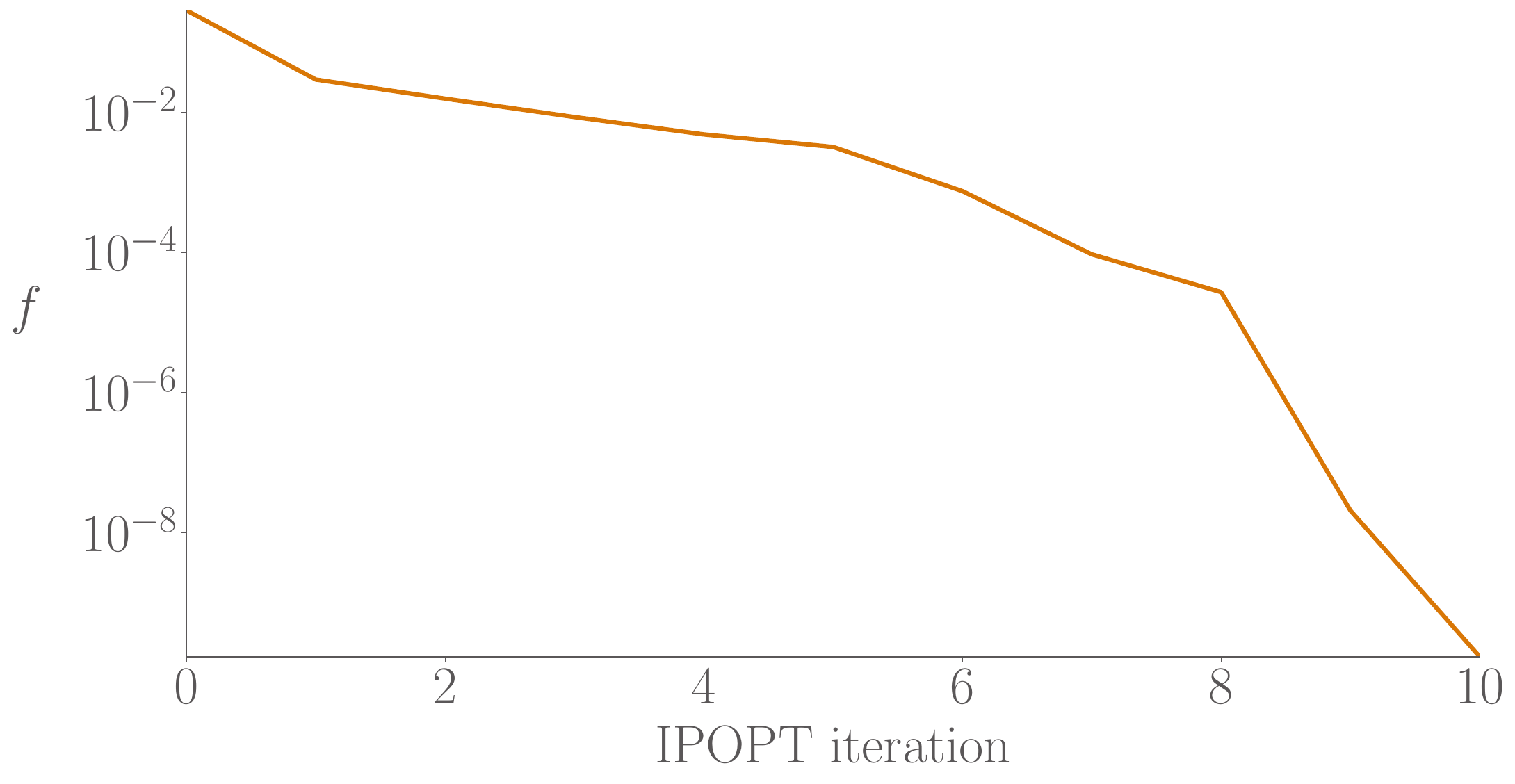}
\caption{Semilog plot of the IPOPT objective history for the objective function $f_1$. 
The adjoint-driven iterations reduce the mode loss by more than eleven orders of magnitude before reaching the prescribed iteration cap.}
\label{fig:burgers1d_obj_history}
\end{figure}

The same adjoint machinery feeds the IPOPT optimizer~\cite{Wchter2005} through pyOptSparse~\cite{Wu2020a} for the mode-matching optimization.
Starting from the uniform control $\boldsymbol{\alpha}_c=\boldsymbol{0}$, the solver converges in ten iterations to $\boldsymbol{\alpha}_c^\star$, with components $(0.25, -0.15, 0.05, 0.15)$ that drive the mode loss from $f_\mathrm{init}=2.83\times10^{-1}$ to $f_\mathrm{opt}=1.70\times10^{-10}$.
Figures~\ref{fig:burgers1d_alpha_cmp}--\ref{fig:burgers1d_obj_history} visualize the Gaussian control profile, the leading POD mode, and the monotone reduction of the objective.
Specifically, Fig.~\ref{fig:burgers1d_mode_cmp} overlays the initial, optimized, and target POD modes and shows that the optimized curve collapses onto the target everywhere except for numerical roundoff; the narrow difference between the initial and optimized curves highlights how the adjoint-based updates steer the modal content.\looseness=-1

The optimized profile aligns with the target Gaussian control distribution extracted from the reference run, and the corresponding mode overlays the target to within machine precision, demonstrating that the adjoint gradient is sufficiently accurate for large-step updates.

In practice, we never know the optimum design conditions.
This is the motivation for optimization.
However, we have manufactured the problem here, and hence know the optimal design vector $\mb{x}=\mb{x}_c^\star=\mbg{\alpha}_c^\star$ before hand.
This enables us to test the robustness of our algorithm from both sides --- the design variable space and state space.

\subsubsection{Different objective functions and multi-mode MCFI for 1D Burger's equation}
\begin{figure}[H]
\centering
\includegraphics[width=1.0\linewidth]{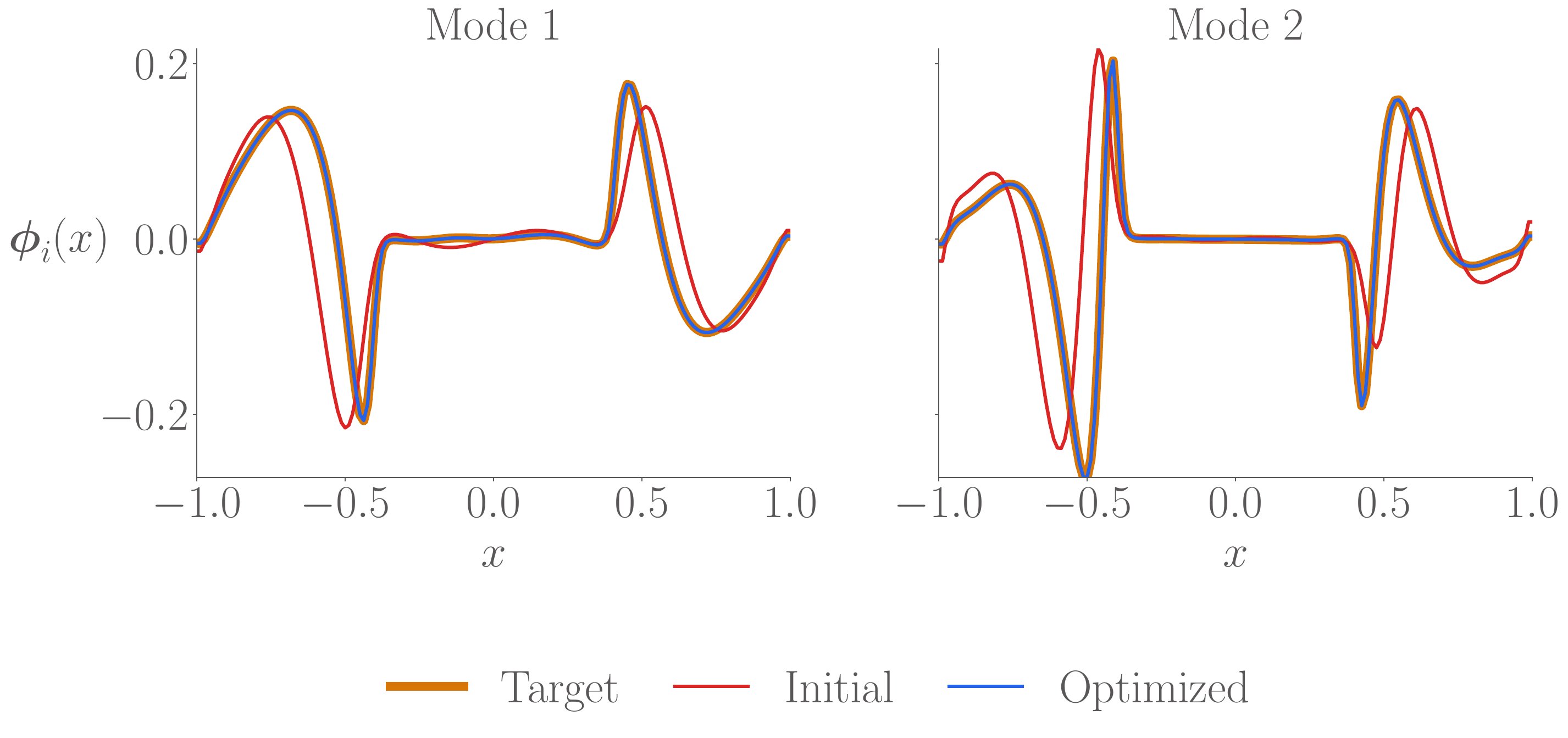}
\caption{Multi-mode comparison for the two-mode objective $f_4$: target modes (thick lines), initial modes, and optimized modes over the spatial domain.
The adjoint formulation applied sequentially to each mode recovers the target shapes to within plotting resolution.}
\label{fig:burgers1d_multi-mode_cmp}
\end{figure}

\begin{figure}[H]
\centering
\includegraphics[width=1.0\linewidth]{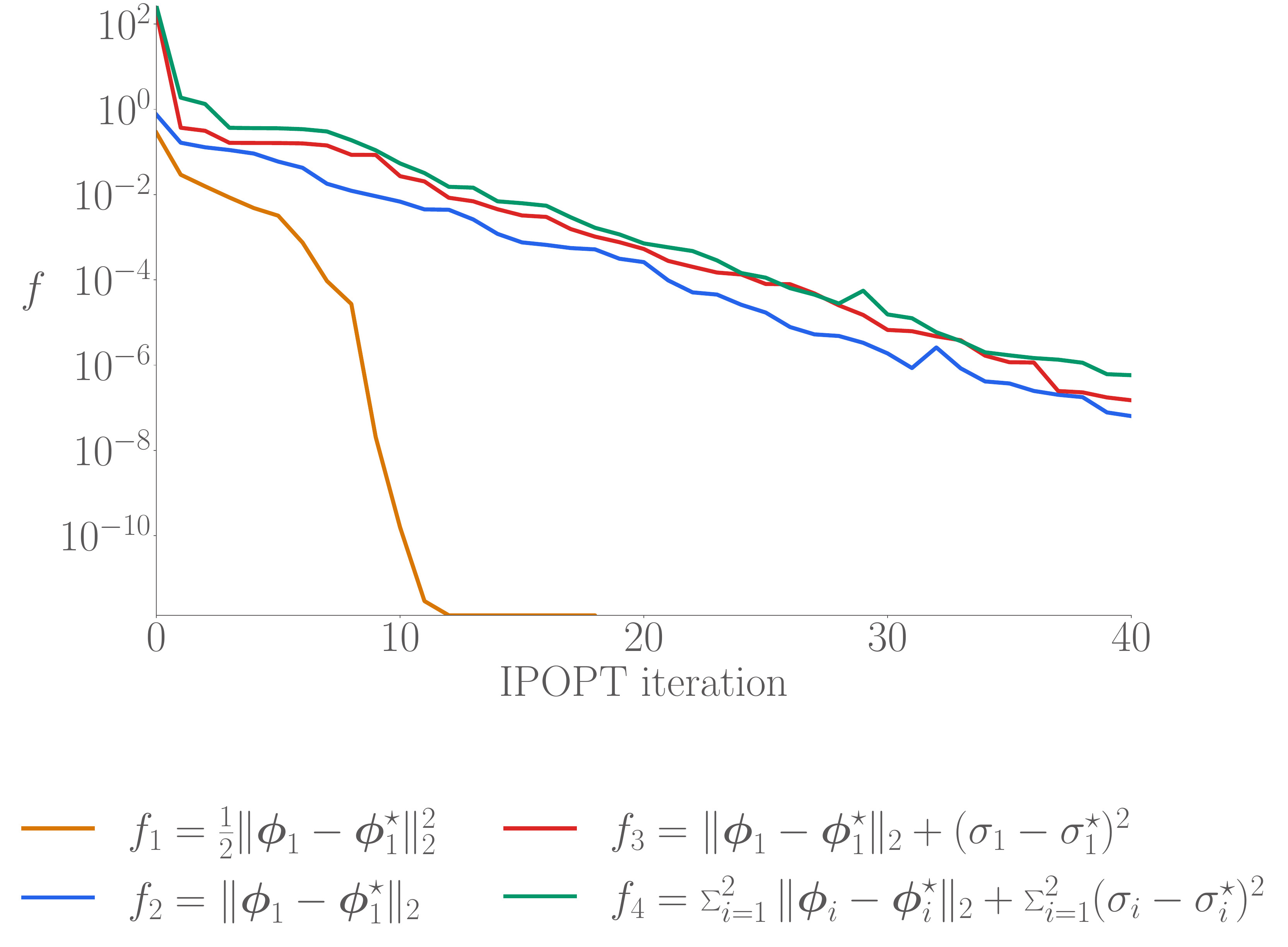}
\caption{Comparison of the IPOPT objective histories for the four functionals defined in Eq.~\Cref{eq:burgers1d_obj_variants}. All cases employ the same adjoint-based sensitivities; the additional POD-singular value penalties in $f_2$, $f_3$ and $f_4$ lengthen the optimization but preserve monotone descent.}
\label{fig:burgers1d_obj_comparison}
\end{figure}

Beyond the objective $f_1$, we explored three additional metrics that extend the mode loss by normalizing, weighting, or augmenting the POD mode loss used during optimization as follows
\begin{equation}\label{eq:burgers1d_obj_variants}
\begin{aligned}
f_1 &= \tfrac{1}{2}\big\|\mbg{\phi}_1 - \mbg{\phi}_1^\star\big\|_2^2, &
f_2 &= \big\|\mbg{\phi}_1 - \mbg{\phi}_1^\star\big\|_2, \\
f_3 &= \big\|\mbg{\phi}_1 - \mbg{\phi}_1^\star\big\|_2 + (\sigma_1 - \sigma_1^\star)^2, &
f_4 &= \sum_{i=1}^{2} \Big(\big\|\mbg{\phi}_i - \mbg{\phi}_i^\star\big\|_2 + (\sigma_i - \sigma_i^\star)^2\Big).
\end{aligned}
\end{equation}
Each function uses the identical adjoint machinery; only the modal forcing assembled from the POD singular triplets changes.
In particular, $f_4$ processes each singular triplet $(\mbg{\phi}_i,\sigma_i,v_i)$ sequentially: the code forms the modal adjoint for mode $i$, evaluates its forcing, and accumulates the dynamical-system adjoint gradients before advancing to the next mode.
The multi-mode outcome for $f_4$ is illustrated in Fig.~\ref{fig:burgers1d_multi-mode_cmp}, which overlays the target, initial, and optimized modes for the two-mode objective.
The optimized modes track the targets closely while preserving smoothness, confirming that the per-mode adjoint solves accumulate consistently in the multi-mode setting.
The IPOPT logs for these four runs are distilled in Fig.~\ref{fig:burgers1d_obj_comparison}, which overlays the iteration histories on a semilog axis.
The purely quadratic loss $f_1$ converges in roughly ten iterations, while the augmented metrics $f_2$, $f_3$ and $f_4$ demand additional iterations to converge.
Nevertheless, all trajectories decrease monotonically, demonstrating that the adjoint gradients remain robust even when the objective couples multiple modal targets.

\subsection{2D Burgers equation}
In this section, we present the results for the forward run for 2D Burgers equation.
Next, we present the verification of the adjoint gradients against those from FD and finally present results for mode-matching optimization.
\subsubsection{Results for the modified unsteady Burgers 2D equation}\label{sec:results2d}
The 2D Burgers PDE was solved numerically with the numerical scheme described in~\Cref{sec:burger_gov}.
The PDE solved is as shown in~\Cref{eq:burgers_u,eq:burgers_v}.
We used the forward Euler time integration as shown in~\Cref{eq:update_compact}.

\begin{figure}[H]
\centering
\includegraphics[width=0.9\linewidth]{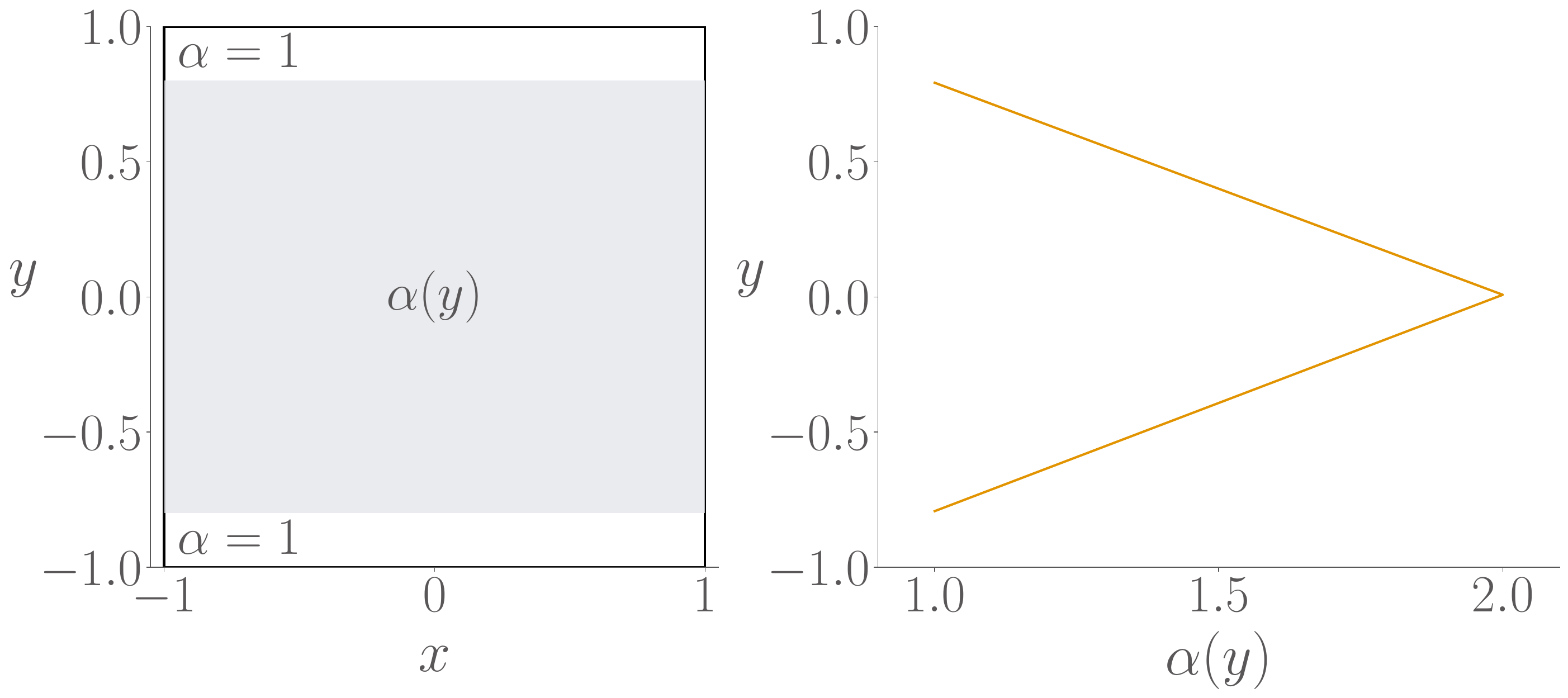}
\caption{Left: schematic of the $[-1,1]^2$ 2D domain used for the 2D Burgers runs. 
Right: strip multiplier $\alpha(y)$.
The profile $\alpha(y)$ shown in the figure on right is the plot for the horizontal strips of constant alpha values in the grey shaded region on the left.
All $\alpha(x,y)$ values in domain outside the grey shaded region are set to a value of 1.}
\label{fig:burgers2d_domain_alpha}
\end{figure}

The computational domain is as shown in~\Cref{fig:burgers2d_domain_alpha}, with ($x$ , $y$) $\in$ $[-1,1]$.
Within $y\in[-0.8,0.8]$ we distribute one hundred horizontal strips and assign each strip a control coefficient $\alpha_{c,j}$.
Outside this band we fix $\alpha=1$.
The ramp profile monotonically increases from unity near $y=-0.8$ to roughly $2$ near $y=0$ before decreasing symmetrically.
These piecewise-constant controls define the advection strength sampled in~\Cref{eq:burgers_vector}.
Thus, we have our design variable vector $\mb{x}$ set equal to the one hundred $\alpha_{c,j}$ values, resulting in one hundred design variables.

The initial condition seeds two opposing Gaussian packets,
\begin{equation}
\label{eq:burgers2d_ic}
\begin{aligned}
u(x,y,0) &= 4\,e^{-\frac{(x-x_L)^2 + (y-y_L)^2}{r_1^2}}
- 4\,e^{-\frac{(x-x_R)^2 + (y-y_R)^2}{r_2^2}}, \\
v(x,y,0) &= 0,
\end{aligned}
\end{equation}
with centers $(x_L,y_L)=(-0.9,0.1)$, $(x_R,y_R)=(0.9,-0.1)$ and radii $r_1=r_2=0.7$.
We normalize $u$ by its maximum absolute value before marching in time.

\begin{figure}[H]
\centering
\includegraphics[width=1.0\linewidth]{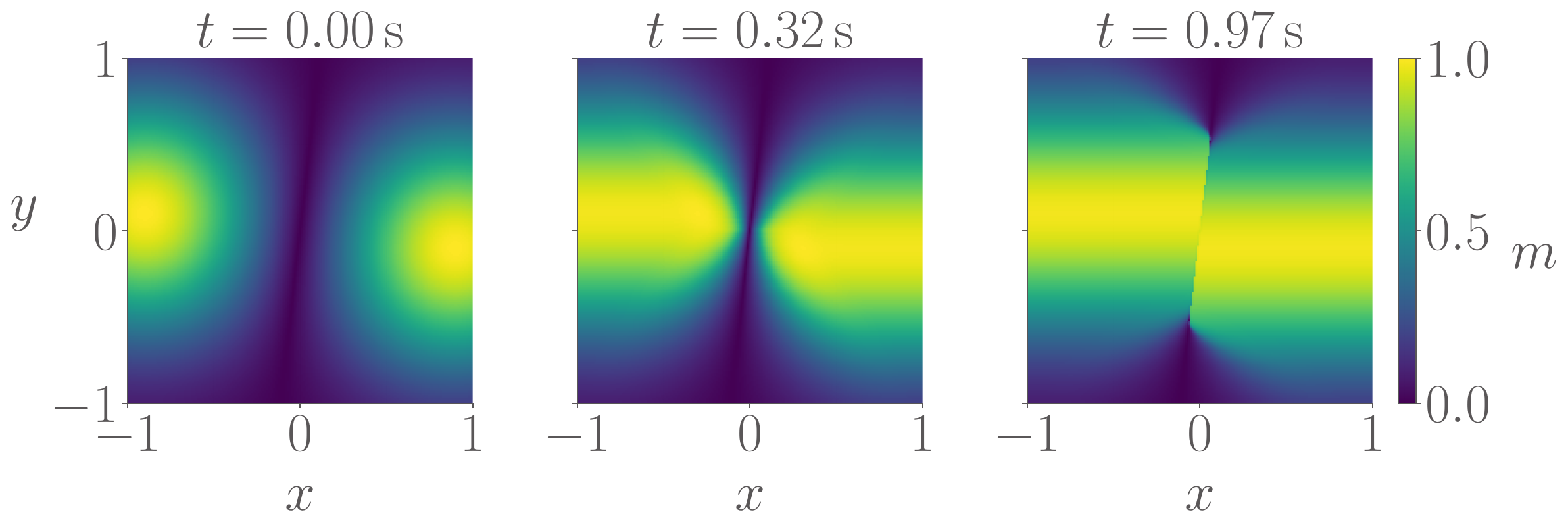}
\caption{Velocity-magnitude ($m=\sqrt{u^2 + v^2}$) snapshots of the baseline 2D Burgers simulation at $t=\{0.00,0.32,0.97\}\,\mathrm{s}$.
The diffused field develops a frontal shock for both Gaussian pulses as they approach each other at prescribed velocities.
The shocks sharpen around 0.32 $\mathrm{s}$.
The shocks interact post 0.32s, leading to a single discontinuous region between roughly $x$ $\in$ $[-0.2, 0.2]$.}
\label{fig:burgers2d_snapshots}
\end{figure}

The domain is discretized with 201 equally spaced points ($n_x=n_y=201$) along $x$ and $y$ directions.
The solver retains the viscosity $\nu=10^{-4}$ and the smooth flux blending from Section~\ref{sec:burger_gov}, so an explicit forward-Euler step is stable under the adaptive CFL restriction
\begin{equation}
\Delta t = \mathrm{CFL}\,\frac{\min(\Delta x,\Delta y)}{\max(|u|)+\max(|v|)+10^{-8}}, \qquad \mathrm{CFL}=0.4,
\end{equation}
with $\Delta x=\Delta y=0.01$.
At the start of the run $\Delta t\approx 4.0\times10^{-3}$ because both velocity components are $O(1)$, and the adaptive step remains in the same range for the remainder of the simulation; reaching $t_{\mathrm{end}}=1.0$ s therefore requires roughly 240 explicit steps.
No manual rescaling of the CFL was needed because the smoothed fluxes keep the Jacobian well conditioned even as the diagonal wave packets interact, and the snapshots highlighted in Fig.~\ref{fig:burgers2d_snapshots} are drawn from this same trajectory.

\subsubsection{Comparison of FD and adjoint method}\label{sec:adj2d}
We applied the adjoint method presented in~\Cref{sec:derivatives}.
The accompanying finite-difference gradients follow the same forward-difference stencil described in~\ref{sec:fd_verification}, with the design vector now representing the strip controls in two dimensions.
The 2D mode-matching problems reuse the same $[-1,1]^2$ domain and horizontal strip band described in~\Cref{sec:results2d}.
\Cref{fig:burgers2d_domain_alpha} summarizes the setup used by the optimizations.
The initial and target $\alpha$ profiles are presented in~\Cref{fig:burgers2d_alpha_profile_obj2} and can be written explicitly as 
\begin{equation}
\alpha_{\text{init}}(y) =
\begin{cases}
1, & |y|>0.8,\\[4pt]
1 + \frac{y+0.8}{0.8}, & -0.8 \le y \le 0,\\[8pt]
1 + \frac{0.8-y}{0.8}, & 0 < y \le 0.8,
\end{cases}
\;\;
\alpha_{\text{tgt}}(y) =
\begin{cases}
1, & |y|>0.8,\\[4pt]
\tilde{\alpha}(y), & 0 \le y \le 0.8,\\[4pt]
\tilde{\alpha}(-y), & -0.8 \le y < 0,
\end{cases}
\label{eq:alpha_profiles_2d}
\end{equation}
where
\begin{equation*}
\tilde{\alpha}(y) =
\begin{cases}
1 + \dfrac{(3.0-1.0)}{0.16}\,y, & 0 \le y < 0.16,\\[6pt]
3.0 - \dfrac{(3.0-0.6)}{0.16}\,(y-0.16), & 0.16 \le y < 0.32,\\[6pt]
0.6 + \dfrac{(2.0-0.6)}{0.16}\,(y-0.32), & 0.32 \le y < 0.48,\\[6pt]
2.0 - \dfrac{(2.0-0.6)}{0.16}\,(y-0.48), & 0.48 \le y < 0.64,\\[6pt]
0.6 + \dfrac{(1.0-0.6)}{0.16}\,(y-0.64), & 0.64 \le y \le 0.8,
\end{cases}
\end{equation*}
after which the piecewise-linear profile is convolved with a discrete Gaussian of width three strips ($\sigma \approx 0.048$) to smooth the strip transitions.

The objective function used for the verification is
\begin{equation}\label{eq:burgers2d_objective}
f=\tfrac{1}{2}\big\|\mbg{\phi}_1(\boldsymbol{\alpha})-\mbg{\phi}_1^\star\big\|_2^2,
\end{equation}
where $\mbg{\phi}_1$ is the first left singular vector of the zero-mean snapshot matrix and $\mbg{\phi}_1^\star$ is is target leading POD mode generated by one forward run with $\boldsymbol{\alpha} = \boldsymbol{\alpha}_\mathrm{tgt}$ from~\Cref{eq:alpha_profiles_2d}.
The adjoint and FD sensitivity evaluations follow exactly the coupled formulation summarized in Sec.~\ref{sec:derivatives} (see especially~\Cref{eq:Adjoint equation,eq:modal_adj_expand,eq:base_adj_expand}.

To ensure consistency with the 1D verification, Table~\ref{tab:2d_grad_check} compares ten randomly selected components of the adjoint gradient against their finite-difference counterparts using the verification run.
It is established that the adjoint-based formulation matches appreciably with the FD gradients.

\begin{table}[H]
\centering
\ra{1.2}
\caption{Adjoint versus finite-difference gradients for ten representative strips in the 2D verification case.}
\label{tab:2d_grad_check}
\vspace{1ex}
\resizebox{\columnwidth}{!}{%
\begin{tabular}{l r r r r}
\toprule
Component &
\multicolumn{1}{c}{Adjoint gradient} &
\multicolumn{1}{c}{FD gradient} &
\multicolumn{1}{c}{$\vert\Delta\vert$} &
\multicolumn{1}{c}{Relative error} \\
\midrule
$\alpha_{c,36}$ & $-0.0022017752901036$ & $-0.0022017752\underline{621117}$ & $2.80\times10^{-11}$ & $1.27\times10^{-8}$ \\
$\alpha_{c,41}$ & $-0.0036555644553539$ & $-0.00365556\underline{50589909}$ & $6.04\times10^{-10}$ & $1.65\times10^{-7}$ \\
$\alpha_{c,63}$ & $-0.0015622453561960$ & $-0.001562245\underline{2914643}$ & $6.47\times10^{-11}$ & $4.14\times10^{-8}$ \\
$\alpha_{c,76}$ & $\phantom{-}0.0040274891413395$ & $\phantom{-}0.004027489\underline{0922287}$ & $4.91\times10^{-11}$ & $1.22\times10^{-8}$ \\
$\alpha_{c,80}$ & $-0.0054657815965023$ & $-0.005465781\underline{6328518}$ & $3.63\times10^{-11}$ & $6.65\times10^{-9}$ \\
$\alpha_{c,83}$ & $-0.0047718753905434$ & $-0.004771875\underline{4642156}$ & $7.37\times10^{-11}$ & $1.54\times10^{-8}$ \\
$\alpha_{c,91}$ & $\phantom{-}0.0035572488683664$ & $\phantom{-}0.003557248\underline{9259783}$ & $5.76\times10^{-11}$ & $1.62\times10^{-8}$ \\
$\alpha_{c,92}$ & $\phantom{-}0.0030171043655690$ & $\phantom{-}0.0030171043\underline{580163}$ & $7.55\times10^{-12}$ & $2.50\times10^{-9}$ \\
$\alpha_{c,93}$ & $\phantom{-}0.0012477984398261$ & $\phantom{-}0.0012477984\underline{062675}$ & $3.36\times10^{-11}$ & $2.69\times10^{-8}$ \\
$\alpha_{c,98}$ & $\phantom{-}0.0001862892115720$ & $\phantom{-}0.000186289\underline{1924054}$ & $1.92\times10^{-11}$ & $1.03\times10^{-7}$ \\
\bottomrule
\end{tabular}%
}
\end{table}

\subsubsection{MCFI for the 2D Burger's equation and optimization}
\begin{figure}[H]
\centering
\includegraphics[width=0.9\linewidth]{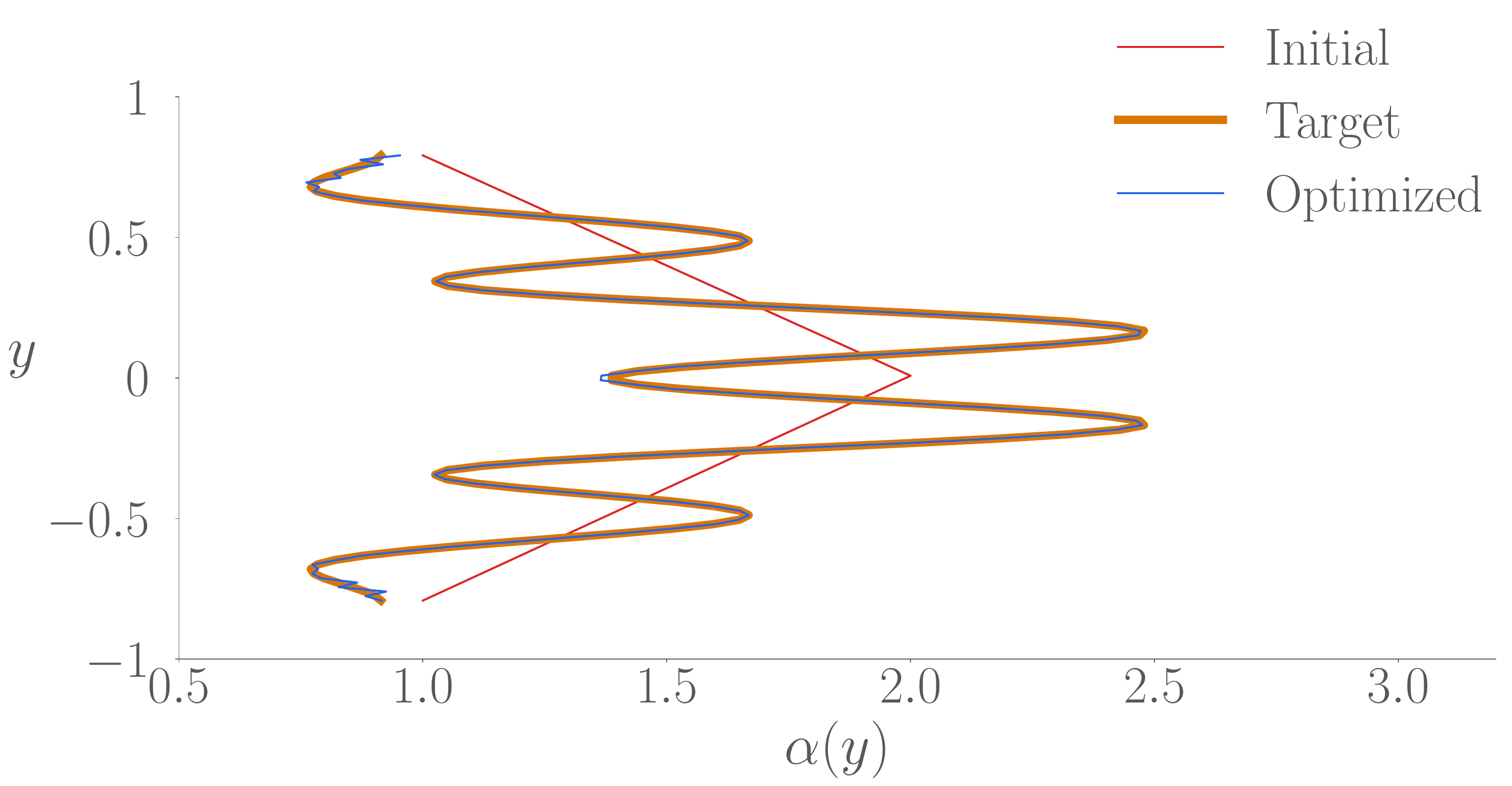}
\caption{Strip-wise $\alpha(y)$ distributions used in the $f_2$ optimization (initial ramp, target multi-peak profile, and optimized design).}
\label{fig:burgers2d_alpha_profile_obj2}
\end{figure}

We investigate two objective functions in this study.
The first is $f_1 = \|\mbg{\phi}_1-\mbg{\phi}^\star_1\|_2$ and second $f_2 = {1}/{2}\|\mbg{\phi}_1-\mbg{\phi}^\star_1\|_2^2$.
The initial, target and optimized $\mbg{\alpha}$ distributions are presented in~\Cref{fig:burgers2d_alpha_profile_obj2} for $f_s$ function.
The results from the two functions are as expected.
Before we jump into the results, we address some observations.

First, it is observed that there is negligible dynamics at play at the extremum $y$ ends of the domain.
This means that the $(\mb{u}\cdot\nabla)\mb{u}$ term in~\Cref{eq:burgers_vector} was roughly zero in that region, leading to minimum influence of $\mbg{\alpha}$ in the region.
Thus, the design variables $\mbg{\alpha}$ in the regions $y \in [-0.8,-0.7]$ and $y \in [0.7,0.8]$ show numerical oscillations and do not match with the target values as the optimizer struggles to find the optimum design parameter there.
The design has almost no influence on the dynamical system state $\mb{U}$, and in turn does not affect the POD mode which depends on the said state variables.

\begin{figure}[H]
\centering
\newcommand{\burgersModeShift}{-1.6cm} 
\includegraphics[width=1.0\linewidth]{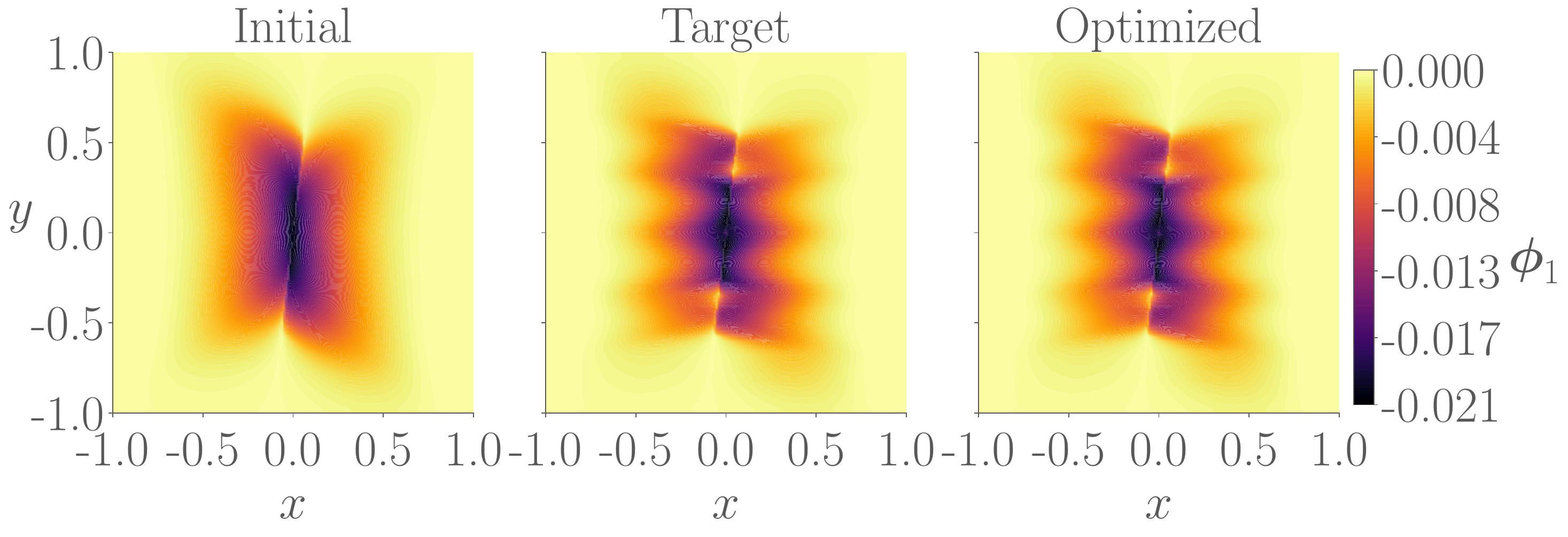}
\caption{Leading POD modes (initial, target, optimized) for the $f_2$ objective. 
The contours are plotted on the physical $x$-$y$ grid using the same smooth color map and scale as the forward solver, making the modal structures comparable to the snapshots in Sec.~\ref{sec:burger_gov}. 
The optimized mode collapses onto the target in the active strip region near $y=0$, while the initial mode retains the bias imposed by the symmetric ramp controller.}
\label{fig:burgers2d_modes_objfunc2}
\end{figure}

Second, in \Cref{fig:burgers2d_modes_objfunc2} however, the optimized mode collapses onto the target everywhere in the domain of influence, while the initial mode retains the bias imposed by the symmetric ramp at $y=0$. 

\begin{figure}[H]
\centering
\includegraphics[width=1.0\linewidth]{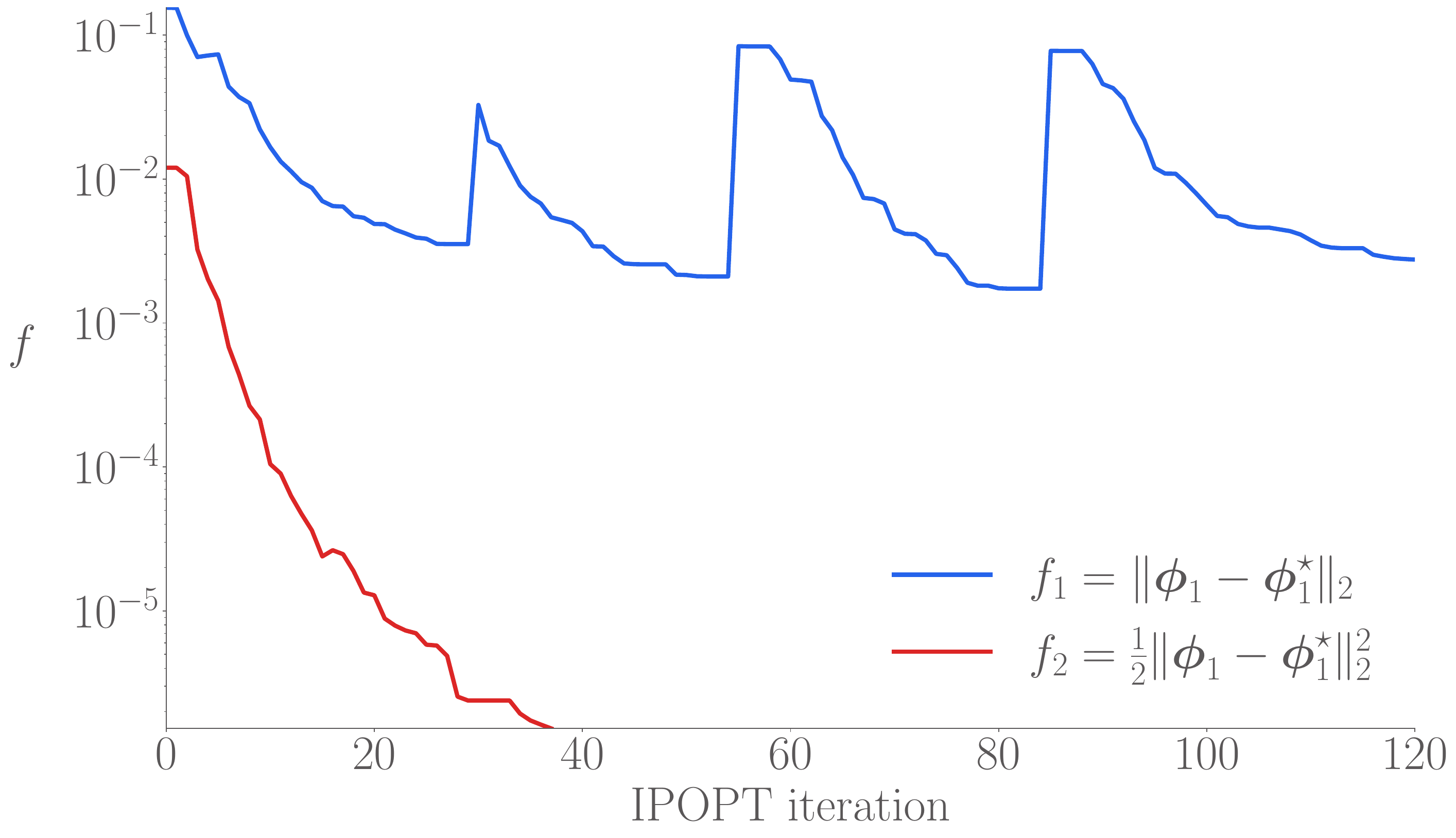}
\caption{Objective histories for the two 2D mode-matching formulations for objective functions $f_1$ and $f_2$.}
\label{fig:burgers2d_objective_history}
\end{figure}

We now discuss the results in conjunction with the aforementioned observations.
The IPOPT histories for $f_1$ and $f_2$ are summarized in \Cref{fig:burgers2d_objective_history}, showing the monotone descent.
The $f_2$ function outperforms the $f_1$ function.
This is a result of the quadratic penalty $f_2 = \tfrac{1}{2}\|\mbg{\phi}_1- \mbg{\phi}_1^\star\|_2^2$ that yields a gradient that scales linearly with the mode loss, providing a smoothly varying curvature that suits line-search globalization.
The gradient of $f_1$ with respect to $\mbg{\phi}_1$ is
\begin{equation}
\f{\p f_1}{\p \mbg{\phi}_1} = \f{\mbg{\phi}_1 - \mbg{\phi}_1^\star}{\|\mbg{\phi}_1-\mbg{\phi}_1^\star\|_2},
\label{eq:df1}
\end{equation} 
and for $f_2$ with respect to $\mbg{\phi}_1$ is
\begin{equation}
\f{\p f_2}{\p \mbg{\phi}_1} = {\mbg{\phi}_1 - \mbg{\phi}_1^\star}.
\label{eq:df2}
\end{equation} 
From~\Cref{eq:df1,eq:df2}, we see two issues.
The first being the sign of the inner product of the current POD mode $\mbg{\phi}_1$ with itself.
This sign must be aligned with that of the target mode.
This is because POD modes are aligned only upto a sign.

The second problem is in~\Cref{eq:df1}.
The denominator in the right hand side of the equation is the L--2 norm of the difference between mode and target mode.
This means that as we approach the target, the denominator tends to zero.

Thus, functions such as $f_1$ are termed $\mathbb{C}^1$ discontinuous and $f_2$ are termed $\mathbb{C}^1$ continuous.
As a result, gradient based optimizers (such as IPOPT~\cite{Wchter2005} used in current study) can take longer steps and reduce the objective more aggressively when minimizing $f_2$--like functions.
Ensuring the sign alignment fixes the only issue with~\Cref{eq:df2}.
Thus, $f_2$ performs better than $f_1$.

Altogether, the results confirm that the adjoint implementation reproduces FD sensitivities, scales to multi-mode objectives, and delivers interpretable optimized controls in both one and two dimensional settings.
With promising results, our future outlook is to extend this to large scale and higher fidelity setups, such as CFD aided optimizations of complex fluid dynamical systems.

\section{Conclusion}\label{sec:conc}

We introduced MCFI, a paradigm that reformulates inverse problems in the reduced space of POD modes rather than the full physical state space.
The core idea—{match structures, not points}—addresses fundamental limitations of conventional field inversion: high dimensionality, ill-conditioning, and lack of inherent regularization when matching full spatio-temporal fields.
By targeting dominant flow structures captured by POD modes, MCFI provides a compact, physically meaningful objective that naturally regularizes the inversion.

Central to this framework is the differentiable POD, an adjoint-based method that efficiently computes sensitivities of POD modes with respect to model parameters.
A key result is that the unsteady adjoint equation needs to be solved only once per simulation, regardless of how many POD modes appear in the objective; modal contributions are simply accumulated into the adjoint right-hand side.
Numerical experiments on 1D and 2D modified Burgers equations confirmed gradient accuracy to near machine precision.

Among the mode-matching objectives evaluated, the quadratic penalty provides the smoothest optimization landscape and fastest convergence while accommodating multiple modes and singular values.
The 2D study demonstrated that optimized spatially varying coefficients concentrate control effort in dynamically active regions, validating MCFI's ability to recover physically meaningful parameter distributions.

These findings establish MCFI as a foundation for scalable inverse design and model calibration in unsteady, high-dimensional systems.
Future work will extend this framework to higher-fidelity CFD solvers and turbulence model calibration.


\bibliographystyle{elsarticle-num-names}
\bibliography{bib/references,bib/extra_bib}

\newpage
\appendix


\section{$(\p \mb{r}_{\text{POD}}/\p \mb{u})^\intercal \mbg{\psi}_\text{POD}$}\label{sec:modal_coupling}
Let $\mb{U} = [\,\mb{u}^{(1)} \;\cdots\; \mb{u}^{(n_t)}\,] \in \mathbb{R}^{n_s\times n_t}$ denote the snapshot matrix assembled from the flattened flow states, where each column $\mb{u}^{(k)} \in \mathbb{R}^{n_s}$ stores one temporal snapshot.
Prior to the POD step we remove the temporal mean from every column.
With $\mb{1} \in \mathbb{R}^{n_t}$ denoting the vector of ones, the centering operator and the zero-mean snapshot matrix are
\begin{equation}
\mb{P} = \mb{I}_{n_t} - \frac{1}{n_t}\,\mb{1}\mb{1}^\intercal,\qquad \tilde{\mb{U}} = \mb{U}\mb{P}.
\label{eq:centering_snap}
\end{equation}
From~\Cref{eq:rmodal}, the blocks that depend on $\tilde{\mb{U}}$ are
\begin{equation}
\mb{r}_{\mathrm{POD},a} = \tilde{\mb{U}}\mb{v} - \sigma\mbg{\phi} \in \mathbb{R}^{n_s},\qquad
\mb{r}_{\mathrm{POD},b} = \tilde{\mb{U}}^\intercal\mbg{\phi} - \sigma\mb{v} \in \mathbb{R}^{n_t},
\end{equation}
while the normalization constraint $\mb{r}_{\mathrm{POD},c} = \mbg{\phi}^\intercal\mbg{\phi} - 1$ is independent of $\tilde{\mb{U}}$.

If we flatten $\tilde{\mb{U}}$ in column-major order, the jacobian of $\mb{r}_{\mathrm{POD}} $ with respect to $\tilde{\mb{U}}$ reads
\begin{equation}
\frac{\partial \mb{r}_{\mathrm{POD}} }{\partial \mathrm{vec}(\tilde{\mb{U}})}
=
\begin{bmatrix}
\mb{v}^\intercal \otimes \mb{I}_{n_s} \\
\mb{I}_{n_t} \otimes \mbg{\phi}^\intercal \\
\mb{0}
\end{bmatrix}
\in \mathbb{R}^{(n_s+n_t+1)\times (n_s n_t)}.
\end{equation}
Let $\boldsymbol{\psi}_{\mathrm{POD}} = [\,\mbg{\psi}_{\mbg{\phi}};\,\mbg{\psi}_{\mb{v}};\,\psi_{\sigma}\,]$ collect the modal adjoint components from~\Cref{eq:modal_adj_expand}. Contracting with $\boldsymbol{\psi}_{\mathrm{POD}}$ gives
\begin{equation}
\left(\frac{\partial \mb{r}_{\mathrm{POD}} }{\partial \mathrm{vec}(\tilde{\mb{U}})}\right)^\intercal
\boldsymbol{\psi}_{\mathrm{POD}}
= \mathrm{vec}\!\left(\mbg{\psi}_{\mbg{\phi}}\mb{v}^\intercal + \mbg{\phi}\mbg{\psi}_{\mb{v}}^\intercal\right),
\end{equation}
because the first block contributes $\mb{v}^\intercal\!\otimes\!\mb{I}_{n_x}$ acting on $\mbg{\psi}_{\mbg{\phi}}$ (yielding $\mbg{\psi}_{\mbg{\phi}}\mb{v}^\intercal$) while the second block contributes $\mb{I}_{n_t}\!\otimes\!\mbg{\phi}^\intercal$ acting on $\mbg{\psi}_{\mb{v}}$ (yielding $\mbg{\phi}\mbg{\psi}_{\mb{v}}^\intercal$). The normalization block is orthogonal to $\psi_{\sigma}$ and does not modify the snapshot coupling. 
Therefore
\begin{equation}
\left(\frac{\partial \mb{r}_{\mathrm{POD}} }{\partial \tilde{\mb{U}}}\right)^\intercal \mbg{\psi}_{\mathrm{POD}}
= \mbg{\psi}_{\mbg{\phi}}\mb{v}^\intercal + \mbg{\phi}\mbg{\psi}_{\mb{v}}^\intercal,
\label{eq:modal_coupling_rhs}
\end{equation}
which supplies the modal forcing term in the dynamical-system adjoint equation.

For multiple modes, the total modal forcing that enters the unsteady adjoint is the sum over modes,
\begin{equation}
\label{eq:multi_modal_forcing}
\mb{g} \;\coloneq\; \sum_{i=1}^{m} \left(\frac{\partial \mb{r}_{\mathrm{POD},i}}{\partial \mathrm{vec}(\tilde{\mb{U}})}\right)^\intercal \mbg{\psi}_{\mathrm{POD},i}
= \sum_{i=1}^{m} \mathrm{vec}\!\left(\mbg{\psi}_{\mbg{\phi},i}\mb{v}_i^\intercal + \mbg{\phi}_i\mbg{\psi}_{\mb{v},i}^\intercal\right),
\end{equation}
where each term follows from~\Cref{eq:modal_coupling_rhs}. This $\mb{g}$ is the quantity accumulated in the algorithmic step that builds the right-hand side of the unsteady adjoint.
Perturbations in $\mb{U}$ propagate through the centering step as $\delta \tilde{\mb{U}} = \delta \mb{U}\,\mb{P}$, which in vector form is $\mathrm{vec}(\delta\tilde{\mb{U}}) = (\mb{P}^\intercal \otimes \mb{I}_{n_x})\,\mathrm{vec}(\delta \mb{U})$.
Applying the chain rule gives the modal residual Jacobian with respect to the original snapshots
\begin{equation}
\frac{\partial \mb{r}_{\mathrm{POD}} }{\partial \mathrm{vec}(\mb{U})}
=
\begin{bmatrix}
(\mb{P}\mb{v})^\intercal \otimes \mb{I}_{n_x} \\
\mb{P}^\intercal \otimes \mbg{\phi}^\intercal \\
\mb{0}
\end{bmatrix}
\in \mathbb{R}^{(n_x+n_t+1)\times (n_x n_t)}.
\label{eq:forcing_term_without_mean}
\end{equation}
In the implementation this projection is achieved by subtracting the temporal mean of $\mbg{\psi}_{\mbg{\phi}}\mb{v}^\intercal + \mbg{\phi}\mbg{\psi}_{\mb{v}}^\intercal$ before forwarding it to the dynamical-system adjoint, which is algebraically equivalent to multiplying by $\mb{P}$ along the time dimension.

\section{Forward FD verification}\label{sec:fd_verification}
To verify the adjoint gradients reported in Secs.~\ref{sec:adj1d} and~\ref{sec:adj2d}, we evaluate a forward FD stencil on the scalar objective $f(\mb{x})$ with respect to the design vector $\mb{x}\in\mathbb{R}^{n_x}$.
Let $\mb{e}_i$ denote the $i^\text{th}$ canonical basis vector and $h_\mathrm{FD}$ the perturbation size.
The FD approximation of the $i^\text{th}$ gradient component reads
\begin{equation}
\left[\nabla_{\mb{x}} f(\mb{x})\right]_i \approx \frac{f(\mb{x} + h_\mathrm{FD}\,\mb{e}_i) - f(\mb{x})}{h_\mathrm{FD}}, \qquad i=1,\ldots,n_x,
\label{eq:fd_forward}
\end{equation}
which is the standard first-order forward-difference formula.
In both the one- and two-dimensional Burgers studies we use $h_\mathrm{FD}=10^{-6}$ and reuse the same function evaluations as the adjoint verification runs.
Because the stencil is written in terms of the abstract design vector $\mb{x}$, it applies unchanged whether $\mb{x}$ represents the four Gaussian-strip coefficients in 1D or the hundred strip controls in 2D.

\section{Residual Jacobian with respect to design variables}\label{sec:residual_jacobian_wrt_x}
For time level $i$ the residual reads
\begin{equation}
\mbg{r}_{\mathrm{uns}}^{(i)} = \mbg{u}^{(i)} - \left[\mbg{u}^{(i-1)} + \Delta t\,\mbg{r}_s\bigl(\mbg{u}^{(i-1)},\mb{x}\bigr)\right],
\end{equation}
so only the spatial residual depends on $\mb{x}$ and
\begin{equation}
\frac{\partial \mbg{r}_{\mathrm{uns}}^{(i)}}{\partial \mb{x}} = -\Delta t\,\frac{\partial \mbg{r}_s(\mbg{u}^{(i-1)},\mb{x})}{\partial \mb{x}} \in \mathbb{R}^{n_s\times n_x},
\end{equation}
where $n_s$ is the number of spatial unknowns in $\mbg{u}$ and $n_x$ is the number of design variables.
Stacking all $n_t$ time levels yields
\begin{equation}
\frac{\partial \mbg{r}_{\mathrm{uns}}}{\partial \mb{x}} =
\begin{bmatrix}
-\Delta t\,\f{\partial\mbg{r}_s(\mbg{u}^{(0)},\mb{x})}{\partial\mb{x}} \\
-\Delta t\,\f{\partial\mbg{r}_s(\mbg{u}^{(1)},\mb{x})}{\partial\mb{x}} \\
\vdots \\
-\Delta t\,\f{\partial\mbg{r}_s(\mbg{u}^{(n_t-1)},\mb{x})}{\partial\mb{x}}
\end{bmatrix}
\in \mathbb{R}^{(n_t n_s)\times n_x}.
\label{eq:res1}
\end{equation}
Each block row is sparse because only the design variables influencing the local stencil produce nonzero columns.
For the modal residual block, we have from~\Cref{eq:rmodal}
\begin{equation}
\frac{\partial \mbg{r}_\mathrm{POD}}{\partial \mb{x}} =
\begin{bmatrix}
\mb{0} \\
\mb{0} \\
\mb{0}
\end{bmatrix}
\in \mathbb{R}^{(2n_s+1)\times n_x},
\label{eq:res2}
\end{equation}
We see a full sparse matrix due to the fact that the modal residuals have no dependence on the design variables in question.
Finally, we perform the global assembly.
The full Jacobian is the vertical concatenation
\begin{equation}
\frac{\partial \mbg{r}_{\mathrm{global}}}{\partial \mb{x}} =
\begin{bmatrix}
\f{\partial \mbg{r}_{\mathrm{uns}}}{\partial \mb{x}} \\
\f{\partial \mbg{r}_\mathrm{POD}}{\partial \mb{x}}
\end{bmatrix}
\in \mathbb{R}^{(n_t n_s + 2n_s + 1)\times n_x}.
\label{eq:res3}
\end{equation}
This matrix multiplies the adjoint vector in the gradient formula $\partial f/\partial \mb{x} - \mbg{\psi}^\intercal (\partial \mbg{r}_{\mathrm{global}}/\partial \mb{x})$ as written in~\Cref{eq:Adjoint equation}, where every column corresponds to one design variable.
This expression gives us the total derivative.
\end{document}